\documentclass[12pt,draft]{amsart}
\usepackage{amsmath,amssymb,amsrefs}

\numberwithin{equation}{section}

\newtheorem{theorem}[equation]{Theorem}
\newtheorem{lemma}[equation]{Lemma}
\newtheorem{corollary}[equation]{Corollary}
\newtheorem{prop}[equation]{Proposition}
\newtheorem{question}[equation]{Question}

\theoremstyle{definition}
\newtheorem{definition}[equation]{Definition}

\theoremstyle{remark}
\newtheorem{remark}[equation]{Remark}

\begin{document}

\begin{center}
\texttt{Comments, suggestions, corrections, and further references
  are most welcomed!}
\end{center}
\bigskip

\title{Capability of some nilpotent products of cyclic groups}
\author{Arturo Magidin}
\address{Department of Mathematical Sciences, The University of
Montana, Missoula MT 59812}
\email{magidin@member.ams.org}

\subjclass[2000]{Primary 20D15, Secondary 20F12}

\begin{abstract}
A group is called capable if it is a central factor
group. We consider the capability of nilpotent products of cyclic
groups, and obtain a generalisation of a theorem of Baer for the small
class case. The approach
may also be used to obtain some recent results on the capability
of certain nilpotent groups of class~$2$. We also obtain a necessary
condition for the capability of an arbitrary $p$-group of class $k$,
and some further results.
\end{abstract}

\maketitle

\section{Introduction}

In his landmark paper on the classification of finite
$p$-groups~\cite{hallpgroups}, P.~Hall remarked that:
\begin{quote}
The question of what conditions a group $G$ must fulfil in order that
it may be the central quotient group of another group $H$, $G\cong
H/Z(H)$, is an interesting one. But while it is easy to write down a
number of necessary conditions, it is not so easy to be sure that they
are sufficient.
\end{quote}
Following M.~Hall and Senior~\cite{hallsenior}, we make the following
definition:

\begin{definition} A group $G$ is said to be \textit{capable} if and only if
there exists a group $H$ such that $G\cong H/Z(H)$; equivalently, if
and only if $G$ is isomorphic to the inner automorphism group of a
group~$H$.
\end{definition}

Capability of groups was first studied by R.~Baer in~\cite{baer},
where, as a corollary of some deeper investigations, he characterised
the capable groups that are a direct sum of cyclic groups.
Capability of groups has received attention in recent
years, thanks to results of Beyl, Felgner, and Schmid~\cite{beyl}
characterising the capability of a group in terms of its epicenter; and
more recently to work of Graham Ellis~\cite{ellis} that describes
the epicenter in terms of the nonabelian tensor square of the
group. This approach was used for example in~\cite{baconkappe} to
characterise the capable $2$-generator finite $p$-groups of odd order
and class~$2$.

While the nonabelian tensor product has proven very useful in the
study of capable groups, it does seem to present certain
limitations. At the end of~\cite{baconkappe}, for example, the authors
note that their methods require ``very explicit knowledge of the
groups'' in question.

In the present work, we will use ``low-tech'' methods to obtain a
number of results on capability of finite $p$-groups. They rely only
on commutator calculus, and so may be a bit more susceptible to
extension than results that require explicit knowledge of the
nonabelian tensor square of a group.

In particular, we will prove a generalisation of Baer's Theorem
characterising the capable direct sums of cyclic groups to
some nilpotent products of cyclic $p$-groups. We will use this result
to derive some consequences and discuss possible directions for
further investigation.

One weakness in our results should be noted explicitly: Baer's Theorem
fully characterises the capable finitely generated abelian groups,
because every finitely generated abelian group can be expressed as a
direct sum of cyclic groups. Our generalisation
does not provide a result of similar
reach, because whenever $k>1$ there exist finite $p$-groups of
class~$k$ that are not a $k$-nilpotent product of cyclic $p$-groups,
even if we restrict to $p>k$.

On the other hand, every finite $p$-group of class $k$ is a
\textit{quotient} of a $k$-nilpotent product of cyclic $p$-groups, so
some progress can be made from our starting point. 

The paper is organized as follows: in Section~\ref{sec:defs} we
present the basic definitions and notation. We proceed in
Section~\ref{sec:nec} to establish a necessary condition for
capability which extends an observation of P.~Hall. We then proceed in
Section~\ref{sec:cent2} to calculate the center of a 2-nilpotent
product of cyclic groups, and derive a new proof of Baer's
Theorem. In
Sections~\ref{sec:cent3} and~\ref{sec:centk} we develop similar
results for the $k$-nilpotent product of cyclic $p$-groups, with
$p\geq k$; from these we obtain in Theorem~\ref{capabilitynilkprod},
the promised generalisation of Baer's Theorem.  In
Section~\ref{sec:kp2} we characterise the capable $2$-nilpotent
products of cyclic $2$-groups. Then in 
Section~\ref{sec:applic} we use our results on capable 2-nilpotent
products of cylic $p$-groups to derive a characterisation of the
capable $2$-generated nilpotent $p$-groups of class two, $p$ an odd
prime (recently obtained through different methods by Bacon
and Kappe), and some other related results, 
by way of illustration of how our method may be
used as a starting point for further investigations. In 
Section~\ref{sec:relquest} we will mention some related questions and
possible future avenues of research. Finally, we conclude in
Appendix~\ref{app:comms} with some results on commutators, which are
used to establish the necessary condition in Section~\ref{sec:nec}. I
have decided to place them to an Appendix because some of the
calculations are involved, and would disrupt the flow of the main
presentation in the paper.

The main results are Theorem~\ref{necessity}, which gives a necessary
condition for capability of a $p$ group of class~$k$, and
Theorem~\ref{capabilitynilkprod}, characterising the capable
$k$-nilpotent products of cyclic $p$-groups for $p>k$.

\section{Definitions and notation\label{sec:defs}}

All maps will be assumed to be group homomorphisms. All groups will be
written multiplicatively, unless we explicitly state otherwise.

Let $G$ be a group. The center of~$G$ is denoted by~$Z(G)$.
The identity element of~$G$ will be denoted by $e$;
if there is some danger of ambiguity, we will use $e_G$.

Let $x\in G$. We say that $x$ is \textit{of exponent
  $n$} if and only if $x^n=e$; we say $x$ is \textit{of order $n>0$}
  if and only if $x^n=e$ and $x^k\neq e$ for all $k$, $0<k<n$. To
  simplify the statements of some of our theorems, we will sometimes
  say that an element is of \textit{order~$0$} to mean that it has
  infinite order.

The commutator $[x,y]$ of two elements $x$ and $y$ is
defined to be $[x,y]=x^{-1}y^{-1}xy$; given two subsets (not necessarily
subgroups) $A$ and~$B$ of~$G$, we let $[A,B]$ be the subgroup
generated by all elements of the form $[a,b]$, with $a\in A$ and $b\in
B$.

The lower central series of~$G$ is the sequence of subgroups defined
by $G_1,G_2,\ldots$, where $G_1=G$, $G_{n+1}=[G_n,G]$. We say a group
$G$ is \textit{nilpotent of class (at most)~$n$} if and only if
$G_{n+1}=\{e\}$; we will often drop the \textit{``at most''} clause, it
being understood. The nilpotent groups of class~$1$ are the abelian
groups.  Note that $G$ is nilpotent of class~$n$ if and only if $G_n\subset
Z(G)$. 

We will write commutators left-normed, so that
\[[a_1,a_2,a_3] = [[a_1,a_2],a_3],\]
etc. The following identities may be verified by direct calculation:
\begin{eqnarray}
\null[xy,z] & = & [x,z][[x,z],y][y,z]\label{prodformone}\\
\null[x,yz] & = & [x,z][z,[y,x]][x,y].\label{prodformtwo}
\end{eqnarray}

The class of all nilpotent groups of class at most~$k$ is denoted
by~$\germ N_k$; it is a variety of groups in the sense of General
Algebra; we refer the reader to Hanna Neumann's excellent
book~\cite{hneumann}. Briefly, a variety of groups is a class of
groups closed under subgroups, quotients, and arbitrary direct
products.  The class of all groups which satisfy the identity $x^n=e$
is denoted $\germ B_n$, and is called the Burnside variety of
exponent~$n$; the class of all nilpotent groups of class at most $k$
and exponent $n$ is ${\germ N_k}\cap {\germ B_n}$.

\subsection*{Nilpotent product of groups}

The nilpotent products of groups were introduced by
Golovin~\cite{golovinnilprods} as examples of regular products of
groups. Although defined in a more general context, in which there is
no restriction on the groups involved, our definition will be
restricted to the situation we are interested in.

\begin{definition} Let $A_1,\ldots,A_n\in{\germ N}_k$. The
\textit{$k$-nilpotent product of $A_1,\ldots,A_n$,} denoted by
$A_1 \amalg^{\germ N_k} \cdots \amalg^{\germ N_k} A_n$,
is defined to be the group $G=F/F_{k+1}$, where $F$ is the free
product of the $A_i$, $F=A_1 * \cdots * A_n$, and $F_{k+1}$ is the
$(k+1)$-st term of the lower central series of $F$.
\end{definition}

Note that the ``$1$-nilpotent product'' is simply the direct
sum. Also, if the $A_i$ are in $\germ N_{k-1}$, and $G$ is the
$k$-nilpotent product of the $A_i$, then the $(k-1)$-nilpotent product
of the $A_i$ is isomorphic to $G/G_k$.

The use of the coproduct notation does not appear to be standard in
the literature, but there is a good reason to use it in our context:

\begin{theorem} Let $A_1,\ldots,A_n \in {\germ N}_k$. The 
  $k$-nilpotent product of the $A_i$ is their \textit{coproduct} in the
  variety $\germ N_k$, in the sense of category theory; that is:
\begin{itemize}
\item[(i)] There exist injective group homomorphism for $j=1,\ldots,n$ 
\[{\rm inc}^{(k)}_j\colon A_j\to A_1\amalg^{\germ
  N_k}\cdots\amalg^{\germ N_k}A_n,\]
 such that the $k$-th nilpotent product is generated
 by the images of the $A_j$.
\item[(ii)] Given any group $K\in{\germ N}_k$, and a family of group
  homomorphisms $\varphi_i\colon A_i\to K$, there exists a unique
  group homomorphism
\[\varphi\colon A_1\amalg^{\germ N_k}\cdots \amalg^{\germ N_k} A_n \to
  K\]
such that for all $i$, $\varphi_i = \varphi\circ {\rm inc}^{(k)}_i$.
\item[(iii)] Condition (i) and (ii) determine $A_1\amalg^{\germ
  N_k}\cdots \amalg^{\germ N_k} A_n$ up to unique isomorphism.
\end{itemize}
\end{theorem}

\begin{proof}
Clause (iii) follows from (i) and (ii) by the usual ``abstract
nonsense'' argument.

For (i), let $i_j\colon A_j\to A_1*\cdots *A_n$ be the canonical
immersion into the free product, and let 
\[\pi\colon A_1*\cdots *A_n \to A_1\amalg^{\germ
  N_k}\cdots\amalg^{\germ N_k}A_n\]
be the canonical projection onto the quotient. Then ${\rm
  inc}^{(k)}_j=\pi\circ i_j$. 

Finally, for (ii), let $\Phi\colon A_1*\cdots *A_n\to K$ be the map
induced by the maps $\varphi_1,\ldots,\varphi_k$, using the universal
property of the free product. Since $K$ is nilpotent of class at most
$k$, the map factors through the $k$-th nilpotent product. It is now
easy to verify that $\varphi=\pi\circ\Phi$ satisfies the given
conditions and is unique.
\end{proof}

Of particular interest is the case when the $A_i$ are cyclic, and
especially when they are also of the same order. Recall that given a
variety of groups ${\germ V}$, a group $F\in {\germ V}$ is 
the \textit{relatively free of rank $n$ in~$\germ V$} if and only if
there exist elements $x_1,\ldots,x_n$ in $F$ such that for any group
$G\in {\germ V}$ and any elements $y_1,\ldots,y_n$ in $G$, there
exists a unique group morphism $\varphi\colon F\to G$ with
$\varphi(x_i)=y_i$. The free groups (in the usual sense) are sometimes
called the \textit{absolutely free groups} for emphasis.

\begin{theorem} Let $A_1,\ldots, A_n$ be cyclic groups, each of order
  $m\geq 0$, generated by $x_1,\ldots,x_n$, respectively.
\item[(1)] If $m=0$, then $A_1\amalg^{\germ N_k}\cdots
  \amalg^{\germ N_k} A_n$ is the relatively free group of rank $n$ in
  $\germ N_k$, freely generated by $x_1,\ldots, x_n$; in particular,
  it is isomorphic to $F(n)/F(n)_{k+1}$, where $F(n)$ is the
  absolutely free group of rank~$n$.
\item[(2)] If $m>0$ and $\gcd(m,k!)=1$, then $A_1\amalg^{\germ
  N_k}\cdots\amalg^{\germ N_k}A_n$ is the relatively free group of
  rank $n$ in ${\germ N_k}\cap{\germ B_m}$, freely generated by
  $x_1,\ldots,x_n$.
\end{theorem}

Condition ``$\gcd(m,k!)=1$'' is equivalent to ``all primes that divide
$m$ are larger than~$k$,'' and it is certainly needed. For example,
the $2$-nilpotent product of two cyclic groups of order $2$ is the
dihedral group of order $8$, and so it cannot be the
relatively free group of rank 2 in ${\germ N_2}\cap{\germ B_2}$,
because the latter variety is just the class of all abelian groups of
exponent $2$.

\subsection*{Basic commutators}

The collection process of M.~Hall gives normal forms for relatively
free groups in ${\germ N_k}$, and in some other special cases. The
concept of basic commutators is essential in this developement. For an
exposition of the collection process, we direct the reader to
Chapter~11 in~\cite{hall}.

\begin{definition} Let $G$ be a group generated by elements
$a_1,\ldots,a_r$. We define the \textit{basic commutators (in
$a_1,\ldots,a_r$)}, their \textit{weights}, and an ordering among
them, as follows:
\begin{itemize}
\item[(1)] $a_1$, $a_2,\ldots,a_r$ are basic commutators of weight
one, and are ordered by the rule $a_1<a_2<\cdots<a_r$.
\item[(2)] If basic commutators of weights less than $n$ have been
  defined and ordered, then $[x,y]$ is a basic commutator of weight
  $n$ if and only if
\begin{itemize}
\item[(a)] $x$ and $y$ are basic commutators, ${\rm wt}(x)+{\rm
  wt}(y)=n$, and
\item[(b)] $x>y$, and
\item[(c)] if $x=[u,v]$, then $y\geq v$.
\end{itemize}
\item[(3)] Commutators of weight $n$ follow all commutators of weight
  less then $n$, and for weight $n$ we define
\[ [x_1,y_1] < [x_2,y_2]\mbox{\ if and only if\ }y_1<y_2\mbox{\ or\
  }y_1=y_2\mbox{\ and\ }x_1<x_2.\]
\end{itemize}
\end{definition}

\begin{theorem}[Basis Theorem; Theorem 11.2.4 in~\cite{hall}]
Let $F$ be the absolutely free group with free generators
$x_1,\ldots,x_r$, and let $c_1,\ldots,c_t$ be the sequence of all
basic commutators in the $x_i$ of weight less than $n+1$, in
nondecreasing order. Then every element $g$ of $F/F_{n+1}$ can be uniquely
expressed as
\[ g = \prod_{i=1}^t c_i^{\alpha_i},\]
where the $\alpha_i$ are integers. Moreover, the basic commutators of
weight~$n$ form a basis for the free abelian group $F_n/F_{n+1}$.
\label{basistheorem}
\end{theorem}

\begin{theorem}[Theorem H1 in \cite{struikone}]
Let $F$, $x_i$, and $c_i$ be as in Theorem~\ref{basistheorem}. If
$g,h\in F/F_{n+1}$ are given by
\[g = \prod_{i=1}^t c_i^{\alpha_i}\qquad\mbox{and}\qquad
h  =  \prod_{i=1}^t c_i^{\beta_i},\]
then
\[ gh = \prod_{i=1}^t c_i^{\gamma_i},\]
where $\gamma_i=f_i(\alpha_j,\beta_k)$ are polynomials with integer
coefficients on the $\alpha_j$ and the $\beta_k$.
\label{struikh1}
\end{theorem}

We also recall a notion related to the weight:

\begin{definition} Let $g\in G$, $g\not= e$. We define $W(g)$ to be
  $W(g)=k$ if and only if $g\in G_k$ and $g\not\in G_{k+1}$. We also
  set $W(e)=\infty$.
\end{definition}

The following properties of the lower central series are
well-known. See for example~\cite{hall}:

\begin{prop} Let $G$ be a group.
\begin{itemize}
\item[(i)] For all $a,b\in G$, $W([a,b])\geq W(a) + W(b)$.
\item[(ii)] If $W(a_i)=w_1$ and $W(b_j)=w_2$, then
\[ \left[ \prod_{i=1}^I a_i^{\alpha_i}\, , \,\prod_{j=1}^J
  b_j^{\beta_j}\right] = \prod_{i=1}^I\prod_{j=1}^J
  [a_i,b_j]^{\alpha_i\beta_j} \pmod {G_{w_1+w_2+1}}.\]
\item[(iii)] If $a\equiv c \pmod{G_{W(a)+1}}$ and $b\equiv d
  \pmod{G_{W(b)+1}}$, then
\[ [a,b] \equiv [c,d] \pmod{G_{W(a)+W(b)+1}}.\]
\item[(iv)] A variant of the Jacobi identity:
\[ [a,b,c]\,[b,c,a]\,[c,a,b] = e \pmod{G_{W(a)+W(b)+W(c)+1}}.\]
\end{itemize}
\label{Widentities}
\end{prop}

Note that if $F$ is the free group on $x_1,\ldots,x_r$, and if $v$ is
a basic commutator on $x_1,\ldots,x_r$ with ${\rm wt}(v)=k$, then
$W(v)=k$ as well. For a general group, we therefore have that ${\rm
wt}(v)\leq W(v)$.

\subsection*{The case of cyclic groups}

In the case of the nilpotent product of cyclic groups, a
result similar to the Basis Theorem can be obtained, at least with
some restrictions on the orders. Explicitly, R.R.~Struik proves:

\begin{theorem}[R.R.~Struik; Theorem 3 in~\cite{struikone}]
Let $A_1,\ldots,A_t$ be cyclic groups of order $\alpha_1,\ldots,\alpha_t$
respectively; if $A_i$ is infinite cyclic, let $\alpha_i=0$. Let $x_i$
generate $A_i$, and let $F$ be the free product
\[ F= A_1 * \cdots *A_t.\]
Let $n\geq 2$ be a fixed positive integer, and assume that all primes
appearing in the factorisations of the $\alpha_i$ are greater than or
equal to $n$. Let $c_1,c_2,\ldots$ be the sequence of basic
commutators on $x_1,\ldots,x_t$ of weight at most $n$. Let
$N_i=\alpha_i$ if $c_i$ is of weight~1, and let
\[ N_i = \gcd(\alpha_{i_1},\ldots,\alpha_{i_k})\]
if $x_{i_j}$, $1\leq j\leq k$, appears in $c_i$. Then every
$g\in F/F_{n+1}$ can be uniquely expressed as
$g = \prod c_i^{\gamma_i}$,
where $\gamma_i$ are integers modulo $N_i$; if $N_i=0$, then
$\gamma_i$ is simply an integer. If
$h = \prod {c_i}^{\delta_i}$
is another element of $F/F_{n+1}$, then
$gh = \prod c_i^{\varepsilon_i}$,
where $\varepsilon_i=f_i(\gamma_j,\delta_k)$ are the polynomials with
integer coefficients of Theorem~\ref{struikh1}.
\label{basisforpcyclic}
\end{theorem}

Thus, at least in the case when the primes involved are sufficiently
large, the basic commutators may be used to give nice normal forms
for the nilpotent products of cyclic groups. 

We also mention the following result:

\begin{theorem}[Ellis, Proposition~6 in~\cite{ellis}]
Let $G=P_1\times \cdots \times P_n$ be a direct product of finitely
generated groups whose abelianizations $P_i^{\rm ab}$ have mutually
coprime exponents. Then $G$ is capable if and only if each $P_i$ is
capable.
\label{reductiontop}
\end{theorem}

Since a torsion nilpotent group is the direct product of its $p$-parts, for
finite nilpotent groups we may use Theorem~\ref{reductiontop} to  
restrict our attention to finite $p$-groups.

Finally, we will also use the following formulas:

\begin{lemma}[Lemma~2 in~\cite{struikone}]
If $G$ is any group and $a,b\in G$, then
\begin{eqnarray*}
\null[a^r,b^s] & \equiv & [a,b]^{rs} [a,b,a]^{s
\binom{r}{2}}
[a,b,b]^{r\binom{s}{2}}
\pmod{G_4}\\
\null[b^r,a^s] &\equiv& [a,b]^{-rs} [a,b,a]^{-r\binom{s}{2}}[a,b,b]^{-s\binom{r}{2}}
\pmod{G_4},
\end{eqnarray*}
where $\binom{r}{2} = \frac{r(r-1)}{2}$ for all integers $r$.
\label{struiklemma2}
\end{lemma}

\subsection*{A note on commutator conventions}

We have defined the commutator with the convention that
$[x,y]=x^{-1}y^{-1}xy$, so that the formula $xy = yx[x,y]$
holds. However, it seems that many if not all recent papers on
capability define it as $[x,y]=xyx^{-1}y^{-1}$, so that $xy = [x,y]yx$
holds. One good reason to use the latter convention would be that it
is the one used for the nonabelian tensor product. The reason we have
not abided by this is just one of simplicity: we will not use the
nonabelian tensor product in this work, and the collection process and
results from~\cite{struikone} we will cite use the former
convention. One could modify slightly the definition of basic
commutators, replace the collection process by ``collecting to the
right,'' write commutators right-normed, and make formulas extend to
the left instead of to the right, and use the latter convention
on the bracket. There is no real problem in doing so, except perhaps
that we are used to thinking of approximation formulas as extending to
the right, not to the left.

\section{A necessary condition\label{sec:nec}}

In this section we give a necessary condition for
capability of finite $p$-groups based on the orders of the elements on a
minimal generating set.  Most of the calculations and hard work
necessary to reach this condition have been placed in
Appendix~\ref{app:comms}. The conclusions of our calculations there
are in Corollary~\ref{likenil2withpowers}.

In the case of small class (that is, when $G\in {\germ N}_k$ is a
$p$-group with $p>k$), the condition reduces to an observation which
goes back at least to P.~Hall (penultimate paragraph in pp.~137
in~\cite{hallpgroups}). Although Hall only considers bases in the
sense of his theory of regular $p$-groups, his argument is essentially
the same as the one we present.  However, Hall's result may not be
very well known, since it is only mentioned in passing; see for
example Theorem~4.4 in~\cite{baconkappe}.

Recall that given a real number~$x$, we let $\lfloor x\rfloor$ denote the
greatest integer smaller than or equal to~$x$.

\begin{lemma} Let $k\geq 1$ be a positive integer, $p$ a prime, and let
  $H$ be a nilpotent $p$-group of class~$k+1$. Suppose that
  $y_1,\ldots,y_r$ are elements of $H$ such that their images generate
  $H/Z(H)$; assume further that the orders of $y_1Z(H),\ldots,y_rZ(H)$
  in $H/Z(H)$ are $p^{\alpha_1},\ldots,p^{\alpha_r}$, respectively,
  with $1\leq \alpha_1\leq\cdots\leq \alpha_r$. Then $\alpha_r\leq
  \alpha_{(r-1)} + \left\lfloor\frac{k-1}{p-1}\right\rfloor$; that is,
\[ \smash{y^{p^{\alpha_{(r-1)}+\left\lfloor\frac{k-1}{p-1}\right\rfloor}}}\in Z(H).\]
\label{acentralelement}
\end{lemma}

\begin{proof}
Since $H$ is generated by $y_1,\ldots,y_r$ and central elements, it
is sufficient to prove that $[y_r^{p^{\alpha_{(r-1)}+\left\lfloor\frac{k-1}{p-1}\right\rfloor}},y_j]=e$ for
$j=1,\ldots,r-1$.  

Note that $y_j^{p^{\alpha_{(r-1)}}}$ is central for $j=1,\ldots,r-1$, so that
\[\forall i\geq \alpha_{(r-1)},\qquad
	  [y_r,y_j^{p^i},y_j]=[y_r,y_j^{p^i},y_r]=e.\]
Thus, we may apply Corollary~\ref{likenil2withpowers} to conclude that
\[ \relax
[y_r^{p^{\alpha_{(r-1)}+\left\lfloor\frac{k-1}{p-1}\right\rfloor}},y_j] =
[y_r,y_j^{p^{\alpha_{(r-1)}+\left\lfloor\frac{k-1}{p-1}\right\rfloor}}]=e,\]
thus proving the lemma.
\end{proof}

The necessary condition is now immediate:

\begin{theorem}[P.~Hall if $k<p$~\cite{hallpgroups}]
Let $k\geq 1$ be a positive integer, and let $p$ be a prime. Let $G$
  be a nilpotent $p$-group of class $k$. Let
  $\{x_1,\ldots,x_r\}$ be a minimal generating set for $G$, and let
  $x_i$ be of order $p^{\alpha_i}$, with $\alpha_1\leq\cdots\leq
  \alpha_r$. If $G$ is capable, then $r>1$ and
\[\alpha_{r}\leq\alpha_{(r-1)}+\left\lfloor\frac{k-1}{p-1}\right\rfloor.\]
\label{necessity}
\end{theorem}

\begin{proof}
Since a center-by-cyclic group is abelian, the necessity of $r>1$ is
clear. So assume that $G$ is capable, and $r>1$. Let $H$ be a
$p$-group of class $k+1$ such that $G\cong H/Z(H)$. Let
$y_1,\ldots,y_r$ be elements of $H$ that project onto $x_1,\ldots,x_r$,
respectively. 
Then Lemma~\ref{acentralelement} gives the condition on $\alpha_r$,
proving the theorem.
\end{proof}

\begin{remark} If $G$ is of small class then 
  $\lfloor\frac{k-1}{p-1}\rfloor=0$, so Theorem~\ref{necessity} says
  that $\alpha_r\leq \alpha_{r-1}$; therefore, when $k<p$ the
  necessary condition becomes ``$r>1$ and
  $\alpha_r=\alpha_{r-1}$.'' This is Hall's observation
  in~\cite{hallpgroups}.
\end{remark}

The necessary condition is reminiscent of Easterfield's bound on the
order of a product. In Theorem~A of~\cite{easterfield}, he proves that
if $x$ and $y$ are any two elements of order $p^{\alpha}$ and
$p^{\beta}$, respectively, of a $p$-groups of class~$c$, then for any
$i\geq \beta+\left\lfloor\frac{c-1}{p-1}\right\rfloor$, we have
$(xy)^{p^i} = x^{p^i}$, and therefore, if we let 
$m=\max\left\{ \alpha,
\beta+\left\lfloor\frac{c-1}{p-1}\right\rfloor\right\}$, then the
order of $xy$ cannot exceed $p^m$. It is
possible that our necessary condition may be derived directly from
Easterfield's work, but I have not attempted to do so: I was unaware
of it until recently.

It would be interesting to know if the inequality is tight.  An easy
case to consider is $p=2$. Then the condition given in
Theorem~\ref{necessity} becomes ``$r>1$ and $\alpha_r\leq
\alpha_{r-1}+(k-1)$.'' In this case, it is easy to find an example 
in which we have equality: the dihedral groups of order
$2^{k+1}$ is of class~$k$, minimally generated by an element of order
$2$ and an element of order $2^k$; and its central quotient is
isomorphic to the dihedral group of order $2^k$. Therefore:

\begin{corollary} For $p=2$, the bound in Theorem~\ref{necessity} is
  tight. That is, for every $k\geq 1$ there exists a capable $2$-group
  of class~$k$,
  with a minimal set of generators that satisfies $\alpha_r=\alpha_{(r-1)}+(k-1)$.
\end{corollary}

Easterfield's paper may also be useful in investigating whether the
given bound is tight; it contains several examples which show many of
the bounds he proves are tight. I hope to investigate this question in
the future.

\section{The center of a $2$-nilpotent product and Baer's
  Theorem\label{sec:cent2}}

In this section we derive a new proof of Baer's characterisation of
the capable of finitely generated abelian groups, by considering the
$2$-nilpotent product of cyclic groups. We present it both because it
provides a rather short proof of Baer's result, and also because it
functions as an introduction to the main ideas that we will use
throughout the rest of the paper.

The case of the $2$-nilpotent products is much simpler than the later
cases we will look at, which allows us to consider more general
situations for the $2$-nilpotent product.
Because of this, we could deal directly with the
original statement of Baer's theorem, which is not restricted to a
\textit{finite} sum of cyclic groups. However, this would introduce
some (easily overcome) notational complications in the proof, which in
turn would lengthen it considerably. Since our main interest
in the sequel is on finite $p$-groups, we will restrict our proof to
the finitely generated case of Baer's Theorem. At the end of the
section we will give the original statement, and indicate how to
modify our proof to obtain it.

The multiplication rules for a $2$-nilpotent product of groups are
straightforward, even without assuming the groups to be cyclic. The
following results are well-known, and we quote them for
information: 

\begin{theorem}[see for example~\cite{machenry}]
 Let $A,B\in\germ{N}_2$. Every element of $A\amalg^{{\germ N}_2} B$ may
 be written uniquely as $\alpha\beta\gamma$, where $\alpha\in A$,
 $\beta\in B$, and $\gamma\in[B,A]$, the cartesian. Moreover, the
 cartesian $[B,A]$ is isomorphic to the tensor product $B^{\rm
 ab}\otimes A^{\rm ab}$ by the map that sends $[b,a]$ to
$\overline{b}\otimes\overline{a}$. 
\end{theorem}

The collection process easily yields:

\begin{theorem} Let $A_1,\ldots,A_r,\ldots\in{\germ N}_2$, and let $G$
 be their $2$-nilpotent product
\[ G = A_1\amalg^{{\germ N}_2} A_2 \amalg^{{\germ N}_2}\cdots
\amalg^{{\germ N}_2} A_r \amalg^{{\germ N}_2}\cdots\]
Every element of $g\in G$ may be written uniquely as
\[ g = \prod_{i=1}^r a_{i} \cdot \prod_{1\leq i<j\leq s}c_{ji},\]
where $r,s$ are nonnegative integers; $a_{i}\in A_{i}$, and
$c_{ji}\in[A_j,A_i]$; for simplicity, we also assume that $s\geq
r$. If $h\in G$ is given by
\[ h = \prod_{i=1}^r b_{i}\cdot \prod_{i\leq i<j\leq
  s}\gamma_{ji},\]
then the product $gh$ is given by
\[ gh = \prod_{i=1}^r (a_{i}b_{i}) \cdot \prod_{1\leq i<j\leq s}
(c_{ji}\gamma_{ji}[a_j,b_i]).\]
\end{theorem}

In the case where the $A_i$ are cyclic groups, the representation
simplifies to yield the description of Theorem~\ref{basisforpcyclic}:

\begin{prop} Let $A_1,\ldots,A_r,\ldots \in {\germ N}_2$ be
  cyclic groups, generated by $x_1,\ldots,x_r,\ldots$
  respectively. Let $\alpha_i$ be the order of $x_i$, and for each
  $i\neq j$, let $\alpha_{ij}=\gcd(\alpha_i,\alpha_j)$. Let $G$ be the
  $2$-nilpotent product of the $A_i$,
\[ G = A_1\amalg^{{\germ N}_2} A_2\amalg^{{\germ N}_2} \cdots
  \amalg^{{\germ N}_2} A_r \amalg^{{\germ N}_2} \cdots\]
Every element $g\in G$ may be written uniquely as
\[ g = \prod_{i=1}^{\infty}x_i^{a_i} \cdot \prod_{1\leq i<j}[x_j,x_i]^{a_{ji}}\]
with $0\leq a_i<\alpha_i$ ($a_i$ an arbitrary integer if $\alpha_i=0$)
and almost all $a_i$ equal to zero; and $0\leq a_{ji}<\alpha_{ji}$
($a_{ji}$ an arbitrary integer if $\alpha_{ji}=0$) and almost all
$a_{ji}$ equal to zero.
\end{prop}

What is the center of the group $G$ described above? 
Let $g\in Z(G)$. Multiplying by suitable elements from $G_2$, we may
assume that $g$ is of the form
$g = x_1^{a_1}\cdots x_r^{a_r}\cdots$, with $0\leq a_i<\alpha_i$
($a_i$ an arbitrary integer if $\alpha_i=0$), and almost all
$a_i=0$. Using the bilinearity of the commutator bracket, we have that
for each $i$:
\[e = [x_i,g]  =  \prod_{1\leq j<i} [x_i,x_j]^{a_j}
\prod_{j>i}[x_j,x_i]^{-a_j}.\] Thus, $g\in Z(G)$ if and only if for
each $i$, for each $j\neq i$ there exists $b_j$ such that
$b_j\equiv a_i \pmod{\alpha_i}$ and $\alpha_{ji}|b_j$. So we obtain:

\begin{theorem}[Baer~\cite{baer}]
Let $G$ be a finitely generated abelian group, and write $G$ as a
direct sum
\[ G = C_1\oplus C_2\oplus\cdots\oplus C_r,\]
where $C_i$ is cyclic of order $\alpha_i$, and
$\alpha_1|\alpha_2|\cdots|\alpha_r$. Then $G$ is capable if and only
if $r>1$ and $\alpha_{r-1}=\alpha_r$.
\end{theorem}

\begin{proof} Necessity follows from Theorem~\ref{necessity}. To prove
  sufficiency, let 
$K = C_1\amalg^{{\germ N}_2} \cdots \amalg^{{\germ N}_2} C_r$, with
  $r>1$ and $\alpha_{r-1}=\alpha_r$. The discussion above
shows that an element of $K$ is
central if and only if it lies in $K_2$, so $K/Z(K)=K/K_2\cong G$.
\end{proof}

The discussion above also proves:

\begin{corollary} Let $K\in {\germ N}_2$, and let
$x_1,x_2,\ldots\in K$ be elements whose orders modulo $Z(K)$ are
$\alpha_1,\alpha_2,\ldots$, respectively. Let
$g = \prod x_i^{a_i}$, and assume that for each $i$, for all $j\neq i$,
there exists an integer $b_j$ such that $b_j\equiv a_i\pmod{\alpha_i}$
and $\alpha_{ij}|b_j$. Then $g\in Z(G)$.
\label{niltwomustbecentral}
\end{corollary}

\begin{remark}
As we noted at the beginning of the section, the original statement of
Baer's Theorem is not restricted to finite direct sums of cyclic
groups. The original result is:
\end{remark}

\noindent\textbf{Baer's Theorem} (Corollary to Existence Theorem
in~\cite{baer})\textbf{.}
\textit{
A direct sum $G$ of cyclic groups (written additively) is capable if
and only if it satisfies the following two conditions:
\begin{itemize}
\item[\textit{(i)}] If the rank of $G/G_{\rm tor}$ is 1, then the orders of the
  elements in $G_{\rm tor}$ are not bounded; and
\item[\textit{(ii)}] If $G=G_{\rm tor}$, and the rank of $(p^{i-1}G)_p/(p^iG)_p$
  is~$1$, then $G$ contains elements of order $p^{i+1}$, for all
  primes~$p$;
\end{itemize}
where $kG=\{kx\mid x\in G\}$, and $H_p = \{h \in H\mid ph=0\}$ for any
subgroup $H$ of~$G$.}

\smallskip
The proof of this result is only complicated by the notation needed to
consider infinite direct sums. When $G/G_{\rm tor}$ is nontrivial, we
write $G$ as a direct sum of infinite cyclic groups and a direct sum
of finite cyclic groups. If there are at least two infinite cyclic
groups, it is easy to verify that the corresponding $2$-nilpotent
product has center equal to its commutator subgroup, and so $G$ will
be capable. If there is exactly one infinite cyclic group, we have two
cases. If every finite cyclic group is of order at most $M>0$,
Corollary~\ref{niltwomustbecentral} implies that if $K$ were a group
with $K/Z(K)\cong G$ then $y^{(M!)}$ would be central, where $y$
projects onto the generator of the infinite cyclic group; so such a
$G$ is not capable. If, on the other hand, the orders of the cyclic
groups are not bounded, the $2$-nilpotent product of the cyclic groups
is a witness to the capability of~$G$. If $G$ is torsion, then it is
convenient to deal with the $p$-parts separately, and so to assume $G$
is a $p$-group. Then we may write
\[ G = \mathop{\oplus}\limits_{j=1}^{\infty}\left(\mathop{\oplus}\limits_{i\in I_j}C_{ij}\right),\]
where $C_{ij}$ is a cyclic group of order $p^{j}$ generated by
$x_{ij}$, and $I_j$ is a (possibly infinite) cardinal. With this
notation, condition (ii) becomes the statement that if $I_j$ is a
singleton, then there exists $k>j$ such that $I_k$ is nonempty.  When
the group satisfies condition (ii) then the corresponding $2$-nilpotent
product of the $C_{ji}$ acts as a witness to the capability of $G$. If
the group does not satisfy condition (ii), let $j_0$ be the last index
with $I_j$ nonempty; then $I_{j_0}$ is a singleton, say $I =
\{0\}$. If $K$ were a group with $K/Z(K)\cong G$, and $x$ projects onto
$x_{0j_0}$, then Corollary~\ref{niltwomustbecentral} proves that
$x_{0j_0}^{p_{j_0}-1}$ is central, yielding a contradiction.

\section{The center of a $3$-nilpotent product\label{sec:cent3}}

In this section we consider the capability of the 2-nilpotent
product of cyclic $p$-groups, $p$ an odd prime, by determining the
center of the $3$-nilpotent product of such cyclic $p$-groups. Once
again, the necessary condition in Theorem~\ref{necessity} is also
sufficient for these 2-nilpotent products.

We start by giving the explicit rules for multiplication in a
$3$-nilpotent product of cyclic groups. From them, we can calculate
the center, which will yield the sufficiency half of
the characterisation theorem.

Let $G$ be the $3$-nilpotent product of cyclic groups,
\[ G = \langle x_1\rangle \amalg^{\germ N_3}\cdots \amalg^{\germ N_3}
\langle x_r\rangle,\]
where $x_i$ is of order $\alpha_i$; assume, moreover, that all
$\alpha_i$ are odd, or equal to $0$ (in the case of infinite cyclic
groups). Let $\alpha_{ij}=\gcd(\alpha_i,\alpha_j)$, and
$\alpha_{ijk}=\gcd(\alpha_i,\alpha_j,\alpha_k)$. We know from
Theorem~\ref{basisforpcyclic} that every element $g$ of $G$ may be written
uniquely as
\begin{eqnarray*} 
g &=& x_1^{a_1}\cdots x_r^{a_r}\prod_{1\leq i<j\leq
  r}[x_j,x_i]^{a_{ji}} \prod_{1\leq i<j\leq
  r}[x_j,x_i,x_i]^{a_{jii}}[x_j,x_i,x_j]^{a_{jij}}\\
&&\qquad\qquad\qquad \prod_{1\leq
  i<j<k\leq r}[x_j,x_i,x_k]^{a_{jik}}[x_k,x_i,x_j]^{a_{kij}}
\end{eqnarray*}
where the $a_i$ are taken modulo $\alpha_i$, the $a_{ji}$ modulo
$\alpha_{ji}$, etc. By considering the case $r=3$ and applying the
collection process, we obtain formulas
  for multiplication of two elements: if $h\in G$ is given by
\begin{eqnarray*}
h &=& x_1^{b_1}\cdots x_r^{b_r}\prod_{1\leq i<j\leq
  r}[x_j,x_i]^{b_{ji}} \prod_{1\leq i<j\leq
  r}[x_j,x_i,x_i]^{b_{jii}}[x_j,x_i,x_j]^{b_{jij}}\\
&&\qquad\qquad\qquad \prod_{1\leq
  i<j<k\leq r}[x_j,x_i,x_k]^{b_{jik}}[x_k,x_i,x_j]^{b_{kij}}
\end{eqnarray*}
then their product $gh$ will be given by
\begin{eqnarray*}
gh &=& x_1^{c_1}\cdots x_r^{c_r}\prod_{1\leq i<j\leq
  r}[x_j,x_i]^{c_{ji}} \prod_{1\leq i<j\leq
  r}[x_j,x_i,x_i]^{c_{jii}}[x_j,x_i,x_j]^{c_{jij}}\\
&&\qquad\qquad\qquad \prod_{1\leq
  i<j<k\leq r}[x_j,x_i,x_k]^{c_{jik}}[x_k,x_i,x_j]^{c_{kij}}
\end{eqnarray*}
where the exponents are given by the following formulas (``${\rm
  mod}\ 0$'' means equality):
\begin{eqnarray*}
c_i & \equiv & a_i+b_i \pmod{\alpha_i};\\
c_{ji} & \equiv & a_{ji} + b_{ji} + a_jb_i \pmod{\alpha_{ji}};\\
c_{jii} & \equiv & a_{jii} + b_{jii} + a_{ji}b_i + a_j\binom{b_i}{2}\pmod{\alpha_{ji}};\\
c_{jij} & \equiv & a_{jij} + b_{jij} + a_{ji}b_j + b_i\binom{a_j}{2} + a_jb_ib_j \pmod{\alpha_{ji}};\\
c_{jik} & \equiv & a_{jik} + b_{jik} + a_{ji}b_k + a_jb_ib_k + a_ja_kb_i -
a_{kj}b_i\pmod{\alpha_{jik}};\\
c_{kij} & \equiv & a_{kij}+b_{kij} + a_{ki}b_j + a_kb_ib_j + a_{kj}b_i\pmod{\alpha_{kij}}.
\end{eqnarray*}
The derivation of these formulas is straightforward, though tedious
and somewhat laborious; they are simply an exercise in the collection
process. We therefore omit it.

These formulas appear in \cite{struikone}, pp.~453. The differences
between those published formulas and the ones here are accounted for
by the fact that Struik uses $[x_i,x_j]$ with $i<j$ rather than its
inverse, and uses $[x_i,x_j,x_k]$ and $[x_j,x_k,x_i]$ instead of
$[x_j,x_i,x_k]$ and $[x_k,x_i,x_j]$ (with $i<j<k$); we can easily
translate from one set of formulas to the other noting that:
\begin{eqnarray*}
\relax[x_j,x_i,x_k] &\equiv& [x_i,x_j,x_k]^{-1} \pmod{G_4}\\
\relax[x_k,x_j,x_i] &\equiv& [x_j,x_i,x_k]^{-1}[x_k,x_i,x_j] \pmod{G_4};
\end{eqnarray*}
both identities follow from Proposition~\ref{Widentities}.

\begin{theorem}
Let $C_1,\ldots,C_r$ be cyclic groups, generated by
$x_1,\ldots,x_r$ respectively. Let $p^{\alpha_i}$ be the order of
$x_i$, where $p$ is an odd prime, and
$1\leq \alpha_1 \leq \alpha_2 \leq \cdots \leq \alpha_r$. Let $G$ be
the $3$-nilpotent product of the $C_i$,
\[ G = C_1 \amalg^{{\germ N}_3} \cdots \amalg^{{\germ N}_3} C_r.\]
Then $Z(G) = \left\langle x^{p^{\alpha_{(r-1)}}}, G_3\right\rangle$.
\label{centerfornil3}
\end{theorem}

\begin{proof} Clearly, $G_3\subset Z(G)$, and $x^{p^{\alpha_{r-1}}}\in
  Z(G)$ follows from Lemma~\ref{acentralelement}.  We need only prove
 the reverse inclusion. Let $\alpha_{ab} =
 \rm{min}(\alpha_a,\alpha_b)$, and analogous with $\alpha_{abc}$.

Suppose that $g\in Z(G)$; since $G_3\subset Z(G)$, we may assume that 
\[ g = x_1^{a_1}\cdots x_r^{a_r} \prod_{1\leq i<j\leq
  r}[x_j,x_i]^{a_{ji}}.\]
First, we compare $gx_r$ with $x_rg$; using the formulas above, we have:
\begin{eqnarray*}
gx_r & = & x_1^{a_1}\cdots x_{r-1}^{a_{r-1}}x_r^{a_r+1} \prod_{1\leq
  i<j\leq r}[x_j,x_i]^{a_{ji}}\\
&& \qquad \prod_{1\leq i<j\leq r}[x_j,x_i,x_r]^{a_{ji}};\\
x_rg & = & x_1^{a_1}\cdots x_{r-1}^{a_{r-1}}x_r^{a_r+1} \prod_{1\leq
  i<j<r} [x_j,x_i]^{a_{ji}} \prod_{1\leq i<r}[x_r,x_i]^{a_{ri}+a_i}\\
&&\qquad \prod_{1\leq i< r}[x_r,x_i,x_i]^{\binom{a_i}{2}}
  [x_r,x_i,x_r]^{a_ia_r}\\
&&\qquad \prod_{1\leq i<j\leq r}[x_r,x_i,x_j]^{a_ia_j}.
\end{eqnarray*}
In particular, by comparing the exponents of $[x_r,x_i]$ in both
expressions we see that we must have $a_i\equiv 0 \pmod{p^{\alpha_i}}$
for $i=1,\ldots,r-1$. Comparing the exponents of $[x_r,x_i,x_j]$ with
$1\leq i<j\leq r$, we see that $a_{ji}\equiv 0
\pmod{p^{\alpha_{ji}}}$. Therefore, we must have $g = x_r^{a_r}$. It
now suffices to show that $p^{\alpha_{r-1}}|a_r$.

To that end, we now calculate $[g,x_{r-1}]$. This must equal $e$,
since $g$ is central, and from Lemma~\ref{struiklemma2} we have:
\begin{eqnarray*}
e & = & [g,x_{r-1}]\\
  & = & [x_r^{a_r},x_{r-1}]\\
  & = & [x_r,x_{r-1}]^{a_r} [x_r,x_{r-1},x_r]^{\binom{a_r}{2}}.
\end{eqnarray*}
Therefore, $a_r\equiv 0
\pmod{p^{\alpha_{r-1}}}$, as desired.
\end{proof}

\begin{corollary}
Let $C_1,\ldots,C_r$ be cyclic $p$-groups, $p$ an odd prime, generated by
$x_1,\ldots,x_r$ respectively. Let $p^{\alpha_i}$ be the order of $x_i$,
and assume that $1\leq \alpha_1 \leq \alpha_2 \leq \cdots \leq \alpha_r$. Let $G$ be
the $2$-nilpotent product of the $C_i$,
$G = C_1 \amalg^{{\germ N}_2} \cdots \amalg^{{\germ N}_2} C_r$.
Then $G$ is capable if and only if $r>1$ and $\alpha_{r-1}=\alpha_r$.
\label{capabilitynil2prod}
\end{corollary}

\begin{proof}
Necessity follows from Theorem~\ref{necessity}. For sufficiency, let
$G$ be given as above with $r>1$ and $\alpha_{r-1}=\alpha_r$. Let $H$
be the $3$-nilpotent product of $C_1,\ldots,C_r$, that is:
$H = C_1\amalg^{{\germ N}_3}\cdots \amalg^{{\germ N}_3} C_r$.
By Theorem~\ref{centerfornil3}, $Z(H)=\langle
H_3,x_r^{\alpha_{r-1}}\rangle = H_3$; therefore, $H/Z(H)$ is simply
the 2-nilpotent product of the $C_i$, so $H/Z(H)\cong G$, as desired.
\end{proof}

\section{The center of a $k$-nilpotent product\label{sec:centk}}

It seems reasonable to guess that the results on the center of a
$k$-nilpotent product of cyclic $p$-groups will hold for arbitrary
$k$, at least as long as $p\geq k$.  Certainly, by developing
multiplication formulas like the ones for the 2- and
3-nilpotent products we could calculate the center explicitly. But that
would involve a rather laborious process.

Luckily, it is possible to get around this complication and employ
instead an inductive argument, which we will do in the present
section.

We begin with two observations on basic commutators.

\begin{lemma}
Let $F$ be the free group on $x_1,\ldots,x_r$. Let $[u,v]$ be a basic
commutator in $x_1,\ldots,x_r$, and assume that ${\rm wt}([u,v])=k$.
\begin{itemize}
\item[(i)] If $v\leq x_r$, then $[u,v,x_r]$ is a basic commutator in
  $x_1,\ldots,x_r$.
\item[(ii)] If $v>x_r$, then
\[ [u,v,x_r] \equiv [v,x_r,u]^{-1}[u,x_r,v]  \pmod{F_{k+2}}.\]
In addition, both $[v,x_r,u]$ and $[u,x_r,v]$ are basic
commutators in $x_1,\ldots,x_r$.
\end{itemize}
\label{jacobigeneral}
\end{lemma}

\begin{proof} Clause (i) follows from the definition of basic
  commutator, as does the claim in Clause (ii) that $[v,x_r,u]$ and
  $[u,x_r,v]$ are basic commutators. The congruence in (ii) follows from
  Proposition~\ref{Widentities}(iv).
\end{proof}

\begin{lemma} Let $F$ be the free group on $x_1,\ldots,x_r$, and let
  $k\geq 2$. Let $c_1,\ldots,c_s$ be the basic commutators in
  $x_1,\ldots,x_r$ of weight exactly $k$ listed in ascending order,
  and write $c_i=[u_i,v_i]$.
  Let $\alpha_1,\ldots,\alpha_s$ be any integers. Let
  $g=c_1^{\alpha_1}\cdots c_s^{\alpha_s}$. For
  $i=1,\ldots,s$, let $d_i$ and $f_i$ be defined by:
\begin{eqnarray*}
d_i &=& \left\{\begin{array}{ll}
\relax[u_i,v_i,x_r]&\mbox{if $v\leq x_r$};\\
\relax[v_i,x_r,u_i]^{-1}&\mbox{if $v>x_r$.}
\end{array}\right.\\
f_i &=&\left\{\begin{array}{cl}
e&\mbox{if $v\leq x_r$}\\
\relax[u_i,x_r,v_i]&\mbox{if $v>x_r$.}
\end{array}\right.
\end{eqnarray*}
Then
\[
\relax [g,x_r]  = \prod_{i=1}^s d_i^{\alpha_i}f_i^{\alpha_i}
\pmod{F_{k+2}}\]
and, except for removing the trivial terms with $f_i=e$,
this expression is in normal form for the abelian group
$F_{k+1}/F_{k+2}$.
\label{wecangodown}
\end{lemma}

\begin{proof} It is easy to verify that if $[u_i,v_i]\neq[u_j,v_j]$,
  then $d_i,f_i,d_j,f_j$ will be pairwise distinct, except perhaps in
  the case where $f_i=f_j=e$. Thus in the expression given for
  $[g,x_r]$, all terms are either powers of the identity or of pairwise
  distinct basic commutators. That the expression does indeed equal
  $[g,x_r]$ follows from Proposition~\ref{Widentities}(ii).
\end{proof}

We are now ready to prove our result:

\begin{theorem} Let $k$ be a positive integer, $p$ a prime, $p\geq k$.
 Let $C_1,\ldots,C_r$
  be cyclic $p$-groups generated
  by $x_1,\ldots,x_r$ respectively. Let $p^{\alpha_i}$ be the order of
  $x_i$, and assume that $1\leq\alpha_1\leq \alpha_2\leq\cdots\leq
  \alpha_r$. If $G$ is
  the $k$-nilpotent product of the $C_i$,
$G = C_1 \amalg^{\germ N_k}\cdots\amalg^{\germ N_k} C_r$,
then $Z(G) = \langle x_r^{p^{\alpha_{(r-1)}}}, G_k\rangle$.
\label{centerknilprod}
\end{theorem}

\begin{proof} One inclusion follows from  Lemma~\ref{acentralelement}
  and the properties of a $k$-nilpotent product. To prove the other
  inclusion, we proceed by induction on $k$. The result is trivially
  true for $k=1$, and we have shown it to be true for $k=2,3$. 

Now assume the result is true for the $(k-1)$-nilpotent product of the
$C_i$, $3<k\leq p$, and let $K=C_1 \amalg^{\germ N_{k-1}}\cdots
\amalg^{\germ N_{k-1}} C_r$; that is, $K=G/G_k$. 

By the induction hypothesis, we know that $Z(K)=\langle
x_r^{p^{\alpha_{(r-1)}}},K_{k-1}\rangle$. Therefore, the center of $G$
is contained in the pullback of this subgroup. So we have the
inclusions
$\langle x_r^{p^{\alpha_{(r-1)}}},G_{k}\rangle \subseteq Z(G)
\subseteq \langle x_r^{p^{\alpha_{(r-1)}}}, G_{k-1}\rangle$.
Let $g\in Z(G)$; multiplying by adequate elements of $G_k$ and an
adequate power of $x_r^{p^{\alpha_{r-1}}}$, we may assume that $g$ is
an element of $G_{k-1}$ which can be written in normal form as:
\[ g = \prod_{i=1}^{n} c_i^{\alpha_i},\]
where $c_{1},\ldots,c_n$ are the basic commutators of weight
exactly $k-1$ in $x_1,\ldots,x_r$, and the $\alpha_i$ are integers on the
relevant interval. If we prove that $g=e$, we will obtain our result.

Since $g\in Z(G)$, its commutator with $x_r$ is trivial. From
Proposition~\ref{Widentities}(ii) we have:
\[e  =  [g,x_r]
  =  \left[ \prod_{i=1}^n c_i^{\alpha_i}, x_r\right]
 =  \prod_{i=1}^{n} [c_i,x_r]^{\alpha_i}.\]

For each $i=1,\ldots,n$, write $c_i=[u_i,v_i]$, with $u_i,v_i$
basic commutators. Let $d_i$ and~$f_i$ be as in the statement of
Lemma~\ref{wecangodown}.  Then we have
\[e = [g,x_r] = \prod_{i=\ell}^n d_i^{\alpha_i}f_i^{\alpha_i}.\]
Except for some $f_i$ which are trivial, the precise ordering of the
remaining terms, and the exponents for the $d_i$ corresponding to
$v_i>x_r$, this is already in normal form. The ordering of the
nontrivial basic commutators is immaterial, since the $d_i$ and $f_j$
commute pairwise; and for those $d_i$ corresponding to $v_i>x_r$, we
simply add the corresponding power of $p$ to the exponent $-\alpha_i$
to obtain an exponent in the correct range.  The expression is then in
normal form, so the only way in which this product can be the trivial
element is if $\alpha_i=0$ for $i=1,\ldots,n$, proving that $g=e$,
and so the theorem.
\end{proof}

From this result we derive immediately, as before:

\begin{theorem} Let $k$ be a positive integer, and let $p$ be a prime greater than~$k$.
Let $C_1,\ldots,C_r$ be cyclic $p$-groups generated by
$x_1,\ldots,x_r$ respectively. Let $p^{\alpha_i}$ be the order of
$x_i$, with $1\leq \alpha_1 \leq \alpha_2 \leq \cdots \leq
\alpha_r$. If $G$ is the $k$-nilpotent product of the $C_i$,
$G = C_1 \amalg^{{\germ N}_k} \cdots \amalg^{{\germ N}_k} C_r$,
then $G$ is capable if and only if $r>1$ and $\alpha_{r-1}=\alpha_r$.
\label{capabilitynilkprod}  
\end{theorem}

\section{The case $k=p=2$\label{sec:kp2}}

In this section we consider the smallest case that is not covered by
our investigations so far: $k=p=2$. In this instance,
Theorem~\ref{necessity} gives the condition $\alpha_r\leq
\alpha_{r-1}+1$. We will prove that the condition is also sufficient
for the case of the $2$-nilpotent product of $2$-groups.

As before, we start by examining the center of a $3$-nilpotent product
of cyclic $2$-groups. Such a product was considered in detail by
R.R.~Struik in~\cites{struikone,struiktwo}, who again obtained both a normal
form and a multiplication table. The main difficulty with simply using
commutator calculus in this situation can be seen by considering the
case when $a$ and $b$ are two elements of order~$2$ in a $3$-nilpotent
group~$G$. Then we have $e = [a,b^2] = [a,b]^2[a,b,b]$, from which we
deduce that $[a,b]^2\in G_3$; this makes it difficult to produce a
normal form using only basic commutators. Struik's solution is to
replace the basic commutators $[z,y,z]$ and $[z,y,y]$ with commutators
$[z^2,y]$ and $[z,y^2]$, respectively, and adjust the ranges of the
exponents accordingly. The normal form result we obtain with these
changes is the following:

\begin{theorem}[Struik; Theorem~4 in~\cite{struikone}]
Let $C_1,\ldots,C_r$ be cyclic groups, generated by $x_1,\ldots,x_r$
respectively, and let the order of $x_i$ be $2^{\alpha_i}$ with $1\leq
\alpha_1\leq
\alpha_2\leq\cdots\leq \alpha_r$. Let $G$ be the $3$-nilpotent product of the
$C_i$,
$G = C_1 \amalg^{{\germ N}_3} C_2 \amalg^{{\germ N}_3}\cdots
\amalg^{{\germ N}_3} C_r$.
Then every element of $g\in G$ can be expressed uniquely in the form
\begin{eqnarray*}
g & = &  x_1^{a_1} x_2^{a_2}\cdots x_r^{a_r} \prod_{1\leq i<j\leq r}
[x_j,x_i]^{a_{ji}} [x_j^2,x_i]^{a_{jij}}[x_j,x_i^{2}]^{a_{jii}}\\
&&\qquad \prod_{1\leq i<j<k\leq r}
[x_j,x_i,x_k]^{a_{jik}}[x_k,x_i,x_j]^{a_{kij}},
\end{eqnarray*}
where the $a_i$, $a_{jik}$ and $a_{kij}$ are integers modulo
$2^{\alpha_i}$; $a_{ji}$ is an integer modulo $2^{\alpha_i+1}$; $a_{jii}$ is an
integer modulo $2^{\alpha_i-1}$, and $a_{jij}$ is an integer modulo
$2^{\alpha_i-1}$ if $\alpha_i=\alpha_j$, and modulo $2^{\alpha_i}$ if $\alpha_i<\alpha_j$. 
\label{nf3nil2group}
\end{theorem}

The multiplication formulas may be obtained either by a direct
calculation with a collection process, or else by using the formulas
for the $3$-nilpotent product given before after suitable manipulation
(see Section~2 in~\cite{struiktwo}). If $g\in G$ is given in 
normal form by:
\begin{eqnarray*}
g & = &  x_1^{c_1} x_2^{c_2}\cdots x_r^{c_r} \prod_{1\leq i<j\leq r}
[x_j,x_i]^{c_{ji}} [x_j^2,x_i]^{c_{jij}}[x_j,x_i^{2}]^{c_{jii}}\\
&&\qquad \prod_{1\leq i<j<k\leq r}
[x_j,x_i,x_k]^{c_{jik}}[x_k,x_i,x_j]^{c_{kij}},
\end{eqnarray*}
and $h\in G$ has normal form:
\begin{eqnarray*}
h & = &  x_1^{d_1} x_2^{d_2}\cdots x_r^{d_r} \prod_{1\leq i<j\leq r}
[x_j,x_i]^{d_{ji}} [x_j^2,x_i]^{d_{jij}}[x_j,x_i^{2}]^{d_{jii}}\\
&&\qquad \prod_{1\leq i<j<k\leq r}
[x_j,x_i,x_k]^{d_{jik}}[x_k,x_i,x_j]^{d_{kij}},
\end{eqnarray*}
then their product is given by
\begin{eqnarray*}
gh & = &  x_1^{f_1} x_2^{f_2}\cdots x_r^{f_r} \prod_{1\leq i<j\leq r}
[x_j,x_i]^{f_{ji}} [x_j^2,x_i]^{f_{jij}}[x_j,x_i^{2}]^{f_{jii}}\\
&&\qquad \prod_{1\leq i<j<k\leq r}
[x_j,x_i,x_k]^{f_{jik}}[x_k,x_i,x_j]^{f_{kij}},
\end{eqnarray*}
where
\begin{eqnarray*}
f_i & = & c_i+d_i;\\
f_{ji} & = & c_{ji}+d_{ji} + c_jd_i - 2\alpha(c_{ji})d_i -
2\alpha(c_{ji})d_j \\
&&\quad -2c_j\binom{d_i}{2} - 2d_i\binom{c_j}{2} - 2c_jd_id_j;\\
f_{jii} & = & c_{jii} + d_{jii} + \alpha(c_{ji})d_i +
c_j\binom{d_i}{2};\\
f_{jij}&=& c_{jij} + d_{jij} + \alpha(c_{ji})d_j + c_jd_id_j +
d_i\binom{c_j}{2};\\
f_{jik} & = & c_{jik} + d_{jik} + \alpha(c_{ji})d_k + c_jd_id_k +
c_jc_kd_i - \alpha(c_{jk})d_i;\\
f_{kij} & = & c_{kij} + d_{kij} + \alpha(c_{ki})d_j + c_kd_id_j +
\alpha(c_{kj})d_i;
\end{eqnarray*}
with $\alpha(c_{ji}) = c_{ji} + 2c_{jii} + 2c_{jij}$; the $f_i$ are
taken modulo appropriate powers of $2$ to place them in the correct range.

Any apparent ambiguity resulting from the different moduli can be seen
to be immaterial. For example, $d_i$ is an integer modulo $2^{\alpha_i}$,
and appears in the formula for $f_{ji}$, which is an integer modulo
$2^{\alpha_i+1}$. But if we replace $d_i$ with $d_i+2^{\alpha_i}$, then the only
possibly ambiguous part of the expression for $f_{ji}$ is
$c_jd_i - 2c_j\binom{d_i}{2}$, which becomes:
\begin{eqnarray*}
\lefteqn{c_j(d_i+2^{\alpha_i}) -2c_j\binom{d_i+2^{\alpha_i}}{2}}\\
&\equiv& c_jd_i + 2^{\alpha_i}c_j - c_j(d_i +
2^{\alpha_i})(d_i+2^{\alpha_i}-1)\pmod{2^{\alpha_i+1}}\\
&\equiv& c_jd_i -c_j(d_i^2-d_i)\pmod{2^{\alpha_i+1}}\\
&\equiv& c_jd_i - 2c_j\binom{d_i}{2}\pmod{2^{\alpha_i+1}},
\end{eqnarray*}
so there is no ambiguity. Once again, our formulas agree with those
labeled (29) in pp.~460 of~\cite{struikone}, once the suitable
modifications are made to account for our use of slightly different
commutators from those used by Struik.

Let $G=C_1\amalg^{{\germ N}_3} C_2 \amalg^{{\germ N}_3}\cdots
\amalg^{{\germ N}_3} C_r$ be the $3$-nilpotent product of cyclic
$2$-groups, with $C_i$ of order $2^{\alpha_i}$, generated by~$x_i$;
we assume that $1\leq \alpha_1\leq\cdots\leq
\alpha_r$. We want to determine the center of such a group. Trivially,
$G_3\subset Z(G)$; so let $g\in Z(G)$. Multiplying by suitable powers
of the commutators of the form $[x_j,x_i,x_k]$ and $[x_k,x_i,x_j]$,
with $1\leq i<j<k\leq r$, we may assume that
\[ g = x_1^{c_1} \cdots x_r^{c_r}\prod_{1\leq i<j\leq r}
   [x_j,x_i]^{c_{ji}} [x_j^2,x_i]^{c_{jij}} [x_j,x_i^2]^{c_{jii}}.\]
We multiply by $h=x_r$ both on the right and the left. 

For $1\leq i<r$, the value of $f_{ri}$ in $gx_r$ is given by
\[
c_{ri}-2\alpha(c_{ri}) \mod 2^{\alpha_i+1},\label{riingh}
\]
while in $x_rg$ it is given by
\[
c_{ri} + c_i - 2\binom{c_i}{2} - 2c_ic_r \mod 2^{\alpha_i+1}.\label{riinhg}
\]
The value of $f_{rii}$ in $gx_r$ is given by
$c_{rii}\mod 2^{\alpha_i-1}$, while in $x_rg$ it is given by
$c_{rii}+\binom{c_i}{2} \mod 2^{\alpha_i-1}$. The value of $f_{rir}$ 
in $gx_r$ is given by 
\[
c_{rir}+\alpha(c_{ri}) \mod\left\{\begin{array}{ll}
2^{\alpha_i}&\mbox{if $\alpha_i<\alpha_r$;}\\
2^{\alpha_i-1}&\mbox{if $\alpha_i=\alpha_r$.}
\end{array}\right.
\]
while the value of $f_{rir}$ in $x_rg$ is:
\[
c_{rir} + c_ic_r \mod\left\{\begin{array}{ll}
2^{\alpha_i}&\mbox{if $\alpha_i<\alpha_r$;}\\
2^{\alpha_i-1}&\mbox{if $\alpha_i=\alpha_r$.}
\end{array}\right.
\]
For $1\leq i<j<r$, the value of $f_{jir}$ in $gx_r$ is given by
$c_{jir} + \alpha(c_{ji})\mod 2^{\alpha_i}$, while the same exponent
in $x_rg$ has value $c_{jir} \mod 2^{\alpha_i}$. 

Finally, in $gx_r$ we have that $f_{rij}$ is equal to 
$c_{rij}\mod 2^{\alpha_i}$, while in $x_rg$ it is equal to
$c_{rij}+c_ic_j\mod 2^{\alpha_i}$.

From these calculations, we obtain the following system of
congruences:
\begin{eqnarray}
-2\alpha(c_{ri}) &\equiv& c_i - 2\binom{c_i}{2}-2c_ic_r
 \pmod{2^{\alpha_i+1}}\label{congone}\\
\binom{c_i}{2}&\equiv& 0 \pmod{2^{\alpha_i-1}}\label{congtwo}\\
\alpha(c_{ri})&\equiv& c_ic_r \quad\left\{\begin{array}{ll}
\left({\rm mod}\ {2^{\alpha_i}}\right)&\mbox{if $\alpha_i<\alpha_r$;}\\
\left({\rm mod}\ {2^{\alpha_i-1}}\right)&\mbox{if $\alpha_i=\alpha_r$.}
\end{array}\right.\label{congthree}\\
\alpha(c_{ji}) & \equiv & 0 \pmod{2^{\alpha_i}}.\label{congfour}\\
c_ic_j&\equiv& 0 \pmod{2^{\alpha_i}}.\label{congfive}
\end{eqnarray}
The first three congruences hold for $i=1,\ldots,r-1$; the last two
hold for each pair $i,j$ satisfying $1\leq i<j<r$. 

From (\ref{congone}) we see $c_i$ is even for $i=1,\ldots,r-1$. From
(\ref{congtwo}) we have $c_i(c_i-1)\equiv 0\pmod{2^{\alpha_i}}$, and since
$0\leq c_i<2^{\alpha_i}$, the only possibility is $c_i=0$; this of
course implies (\ref{congfive}). Combining this
with (\ref{congone}) we conclude that $\alpha(c_{ri})\equiv 0
\pmod{2^{\alpha_i}}$ for $i=1,\ldots,r-1$; together with 
(\ref{congfour}) this gives that $\alpha(c_{ji})\equiv 0
\pmod{2^{\alpha_i}}$ for all $1\leq i<j\leq r$.
Therefore, $g$ can be written as
\[ g = x_r^{c_r} \prod_{1\leq i<j\leq r}
   [x_j,x_i]^{c_{ji}}[x_j,x_i^2]^{c_{jii}}[x_j^2,x_i]^{c_{jij}},\]
subject to the condition that 
\[ \alpha(c_{ji}) = c_{ji} + 2c_{jii} + 2c_{jij} \equiv 0
\pmod{2^{\alpha_i}},\quad1\leq i<j\leq r.\]

Next we look at the result of multiplying $g$ by $[x_r,x_{r-1}]$; in
$g[x_r,x_{r-1}]$, we have $f_{r,(r-1)}=c_{r,(r-1)}+1 \mod
2^{\alpha_{(r-1)}+1}$; while in $[x_r,x_{r-1}]g$ we have
$f_{r,(r-1)}=c_{r,(r-1)}+1-c_r\mod 2^{\alpha_{(r-1)}+1}$. Since $g$ is
central, we conclude that $c_r\equiv 0 \pmod{2^{\alpha_{(r-1)}+1}}$.

These conditions are both necessary and sufficient. We obtain:

\begin{theorem}
Let $C_1,\ldots,C_r$ be cyclic groups, generated by $x_1,\ldots,x_r$
respectively; assume that the order of $x_i$ is $2^{\alpha_i}$, and
that the exponents satisfy 
$1\leq \alpha_1\leq \alpha_2\leq \cdots \leq \alpha_r$. Let $G$ be the
3-nilpotent product of the $C_i$,
$G = C_1\amalg^{{\germ N}_3}\cdots \amalg^{{\germ N}_3} C_r$.
Then $g\in Z(G)$ if and only if it can be written in normal form as:
\begin{eqnarray*}
 g &=& x_r^{c_r} \prod_{1\leq i<j\leq r}
 [x_j,x_i]^{c_{ji}}[x_j,x_i^2]^{c_{jii}}[x_j^2,x_i]^{c_{jij}}\\
&&\quad \prod_{1\leq i<j<k\leq r}
 [x_j,x_i,x_k]^{c_{jik}}[x_k,x_i,x_j]^{c_{kij}},
\end{eqnarray*}
where $c_r\equiv 0 \pmod{2^{\alpha_{(r-1)}+1}}$, and
$\alpha(c_{ji})\equiv 0 \pmod{2^{\alpha_i}}$. That is, $Z(G)$ is
generated by $x_r^{2^{\alpha_{(r-1)}+1}}$; for each $1\leq i<j\leq r$ the elements 
\[ [x_j,x_i]^{2^{\alpha_i+1}-2}[x_j,x_i^2],\quad
   [x_j,x_i]^{2^{\alpha_i+1}-2}[x_j^2,x_i],\quad\mbox{and}\quad [x_j,x_i]^{2^{\alpha_i}};\]
and the elements $[x_j,x_i,x_k]$ and $[x_k,x_i,x_j]$ with $1\leq
   i<j<k\leq r$. In other words,
\[ Z(G) = \Biggl\langle \Bigl\{x_r^{2^{\alpha_{(r-1)}+1}},
   G_3\Bigr\}\;\bigcup\; \Bigl\{
   [x_j,x_i]^{2^{\alpha_i}}\,\Bigm|\, 1\leq i<j\leq r\Bigr\}\Biggr\rangle.\]
\label{centerofnil32groups}
\end{theorem}

\begin{proof}
The necessity comes from our previous discussion. 
Sufficiency is simply a matter of verifying that for $i=1,\ldots,r$,
$gx_i=x_ig$, which is straightforward and will be ommitted. As for our
choice of generators, they are certainly all in the center as
described; we have chosen them because, for example,
$\relax[x_j,x_i]^{2^{\alpha_i+1}-2}[x_j,x_i^2] = [x_j,x_i,x_i]$.
It is now easy to verify every element satisfying the conditions given may
be expressed as a product of the given elements. 
\end{proof}

With a description of the center, we can now easily derive the
characterisation of the capable $2$-nilpotent products of cyclic
$2$-groups: 

\begin{theorem}
Let $C_1,\ldots,C_r$ be cyclic $2$-groups, and let
$x_1,\ldots,x_r$ be their respective generators. Let $2^{\alpha_i}$ be the order of
$x_i$, and assume that $1\leq \alpha_1\leq\cdots\leq \alpha_r$. If
$G$ is the $2$-nilpotent product of the $C_i$,
\[ G = C_1 \amalg^{{\germ N}_2} C_2 \amalg^{{\germ N}_2} \cdots
\amalg^{{\germ N}_2} C_r,\]
then $G$ is capable if and only if $r>1$ and $\alpha_r \leq \alpha_{(r-1)}+1$.
\end{theorem}

\begin{proof} Necessity follows from Theorem~\ref{necessity}. For
  sufficiency, let $K$ be the $3$-nilpotent product of the $C_i$,
$K = C_1 \amalg^{{\germ N}_3} C_2 \amalg^{{\germ N}_3} \cdots
\amalg^{{\germ N}_3} C_r$.
Then the description of the center at the end of
Theorem~\ref{centerofnil32groups} makes it easy to verify that
$K/Z(K)\cong G$, so $G$ is capable.
\end{proof} 

\section{Some applications of our approach\label{sec:applic}}

As we noted in the introduction, one weakness of
Theorem~\ref{capabilitynilkprod} is that whereas the $1$-nilpotent
products of cyclic groups covers all finitely generated abelian
groups, the case $k\geq2$ does not do the same for the finitely
generated nilpotent groups of class~$k$. However, it is possible to
use our results as a starting point for discussing capability of other
more general $p$-groups. We present here one example.

A recent result of Bacon and Kappe characterises the capable
2-generated nilpotent $p$-groups of class two with $p$ an odd prime
using the nonabelian tensor square. We can recover their result using
our techniques, and obtain a bit more.

In Theorem~2.4 of~\cite{baconkappenonab}, the authors present a
classification of the finite $2$-generator $p$-groups of class two,
$p$ an odd prime. With a view towards their calculations of the
nonabelian tensor square, the authors classify the groups into three
families. We will modify their classification and coalesce them into a
single presentation.

Let $G=\langle a,b\rangle$ be a finite nonabelian $2$-generator
$p$-group of class~$2$, $p$ an odd prime. Then $G$ is isomorphic to
the group presented by:
\begin{equation}
\left\langle a,b\,\Biggm|\,
\begin{array}{rcl}
a^{p^{\alpha}}=b^{p^{\beta}}=[b,a]^{p^{\gamma}} &=&e,\\
\relax[a,b,a]=[a,b,b] &=& e,\\
a^{p^{\alpha+\sigma-\gamma}} [b,a]^{p^{\sigma}}& = & e.
\end{array} \right\rangle
\label{eq:presentation}
\end{equation}
where $\alpha+\sigma\geq 2\gamma$, $\beta\geq\gamma\geq 1$,
$\alpha\geq\gamma$, and if $\sigma=\gamma$, then $\alpha\geq\beta$.
Under these restrictions, the choice is uniquely determined.

From the above conditions, we get that $0\leq \sigma\leq\gamma$. If
$\sigma=\gamma$, we obtain the groups in Bacon and Kappe's first
family, which one might call the ``coproduct like'' groups (they are
obtained from the nilpotent product $\langle a\rangle\amalg^{{\germ
    N}_2}\langle b\rangle$ by moding out by a power of $[a,b]$). If
$\sigma=0$, we obtain the split meta-cyclic groups, which are the
second family in~\cite{baconkappenonab}. The cases $0<\sigma<\gamma$
correspond to their third family.

Their result, which appears in~\cite{baconkappe}, is that a
$2$-generated group with presentation as in~(\ref{eq:presentation})
with $\sigma=0$ or $\sigma=\gamma$ is capable if and only if
$\alpha=\beta$. The condition is also both necessary and sufficient
for the remaining case with $0<\sigma<\gamma$ (in this
case,~\cite{baconkappe} contains an error which the
authors are in the process of correcting~\cite{kappepers}).

Thus, in the case of $2$-generated $p$-groups of class~two, $p$ an odd
prime, Baer's condition is both necessary and sufficient, just as for
finite abelian groups.

Although we could prove the result with just one argument, we will
divide it in two in order to prove slightly more for the case where
$\sigma=\gamma$. 

Let $C_1,\ldots,C_r$ be cyclic groups generated by $x_1,\ldots,x_r$
respectively. Let $p^{\alpha_i}$ be the order of $x_i$, and assume
that $p$ is an odd prime, and $1\leq\alpha_1\leq\alpha_2\leq\cdots
\leq\alpha_r$.  Let $K$ be the $2$-nilpotent product of the $C_i$,
\[ K = C_1 \amalg^{{\germ N}_2}\cdots \amalg^{{\germ N}_2} C_r.\]
Let $\alpha_{ji}=\alpha_i$ for $1\leq i<j\leq r$, and
for each pair $j,i$, let $\beta_{ji}$ be a positive integer less than
or equal to~$\alpha_{ji}$. Let $N=\langle
[x_j,x_i]^{p^{\beta_{ji}}}\rangle$; since $N$ is central, we have
$N\triangleleft K$. Let $G=K/N$.

\begin{theorem} Notation as in the previous paragraph. $G$ is capable
  if and only if $r>1$ and $\alpha_{r-1}=\alpha_r$.
\end{theorem}

\begin{proof} Necessity follows from Theorem~\ref{necessity}. For
 sufficiency, let $H$ be the $3$-nilpotent product of the $C_i$,
$H = C_1\amalg^{{\germ N}_3}\cdots \amalg^{{\germ N}_3} C_r$.
Clearly $H$ will not do, so we need to take a quotient of $H$ so that,
in the resulting group, $[x_j,x_i]^{\beta_{ji}}$ is central. To that
end, we let $M$ be the subgroup of $H$ generated by all elements of the form
$[x_j,x_i,x_k]^{\beta_{ji}}$ with $1\leq i<j\leq r$, and $k$ arbitrary.
In terms of the normal forms we have given, this means the elements
$[x_j,x_i,x_k]^{\beta_{ji}}$ for $k\geq i$, and the elements
\[[x_i,x_k,x_j]^{-\beta_{ji}}[x_j,x_k,x_i]^{\beta_{ji}}\]
for $k<i$.  It is easy to verify that all elements of~$H/M$ have a
normal form as in Theorem~\ref{basisforpcyclic}, and that the
multiplication of these elements uses the same formulas as those in
$H$, except that the exponents of the basic commutators of weight~$3$
are now taken modulo the adequate $\beta_{ji}$ instead of
$\alpha_{ji}$.

Proceeding now as in the proof of Theorem~\ref{centerfornil3}, one proves
that the center of $H/M$ is generated by the third term
of the lower central series, the image of $x_r^{p^{\alpha_{(r-1)}}}$,
and those of the elements
$[x_j,x_i]^{p^{\beta_{ji}}}$; therefore, moding out by the center, we
obtain the group~$G$, as desired.
\end{proof}

The case where $\beta_{ji}=\alpha_{ji}$ for each $1\leq i<j\leq n$
corresponds to the first part of Corollary~4.3 in~\cite{baconkappe}.
Now we consider the case with $\sigma<\gamma$.

\begin{theorem}[cf. Corollary 4.3 in~\cite{baconkappe}]
 Let $p$ be an odd prime, and let $G$ be a group presented by
 $(\ref{eq:presentation})$, with $0\leq\sigma<\gamma$. Then $G$ is capable
 if and only if $\alpha=\beta$.
\end{theorem}

\begin{proof} Necessity once again follows from Theorem~\ref{necessity}; so we
  only need to prove sufficiency.  We will construct the ``obvious''
witness to the capability, by starting with the $3$-nilpotent product
of two cyclic groups of order $p^{\alpha}$, generated by~$x$ and~$y$;
then for every relation $r$ in the presentation of~$G$, we will mod
out by the subgroup $\langle [r,x],[r,y]\rangle$, thus making making
$r$ central. Then we just need to make sure that in the resulting
group, the map $x\mapsto a$, $y\mapsto b$ will yield the desired
isomorphism between the central quotient and~$G$.

(In essence, what we are doing is constructing a ``generalised 
extension of~$G$''  which can be used to
determine the capability of~$G$. See Theorem~III.3.9 in~\cite{beyltappe})

So let $K_0 = \langle x\rangle \amalg^{{\germ N}_3} \langle y\rangle$,
where $x$ and $y$ are both of order $p^{\alpha}$. 

First, we want to make sure that $[y,x]^{p^{\gamma}}$ is central, so
we let $N=\langle [y,x,x]^{p^{\gamma}},[y,x,y]^{p^{\gamma}}\rangle$,
and let $K_1=K_0/N$.

The next step is to ensure that
$x^{p^{\alpha+\sigma-\gamma}}[y,x]^{p^{\sigma}}$ is central. So first
we consider
\begin{eqnarray*}
[x^{p^{\alpha+\sigma-\gamma}}[y,x]^{p^{\sigma}},x] & = &
[x^{p^{\alpha+\sigma-\gamma}},x] [x^{p^{\alpha+\sigma-\gamma}},x,[y,x]^{p^{\sigma}}]
[[y,x]^{p^{\sigma}},x]\\
& = & [y,x,x]^{p^{\sigma}}
\end{eqnarray*}
(using (\ref{prodformone}) and Proposition~\ref{Widentities}). So we
let $K_2=K_1/\langle [y,x,x]^{p^{\sigma}}$. 

With these quotients, the only difference in the normal form and
multiplication tables for $K_2$ and for $K_0$ is the order of
$[y,x,x]$ and $[y,x,y]$.

The final quotient we need to take is to ensure that
$x^{p^{\alpha+\sigma-\gamma}}[y,x]^{p^{\sigma}}$ also commutes with
$y$. To that end, we consider:
\begin{eqnarray*}
[x^{p^{\alpha+\sigma-\gamma}}[y,x]^{p^{\sigma}},y] & = &
 [x^{p^{\alpha+\sigma-\gamma}},y] [x^{p^{\alpha+\sigma-\gamma}},y,[y,x]^{p^{\sigma}}]
[[y,x]^{p^{\sigma}},y]\\
& = & [x^{p^{\alpha+\sigma-\gamma}},y] [y,x,y]^{p^{\sigma}}\\
& = & [y,x]^{-p^{\alpha+\sigma-\gamma}}
 [y,x,x]^{-\binom{p^{\alpha+\sigma-\gamma}}{2}} [y,x,y]^{p^{\sigma}}.
\end{eqnarray*}
The last equality uses Lemma~\ref{struiklemma2}. Note, however, that
the conditions on $\alpha$, $\sigma$, and~$\gamma$ imply that
$\alpha+\sigma-\gamma>\sigma$, so we can
simplify the expression above to:
\[ [x^{p^{\alpha+\sigma-\gamma}}[y,x]^{p^{\sigma}},y] 
= [y,x]^{-p^{\alpha+\sigma-\gamma}}[y,x,y]^{p^{\sigma}}.\]
So let $N=\langle
[y,x]^{p^{\alpha+\sigma-\gamma}}[y,x,y]^{-p^{\sigma}}\rangle$. 
Using the multiplication formulas for $K_2$, it is easy to verify that
for all $g\in K_2$, $[g,x],[g,y]\in N$ if and only if 
\[  g \in \left\langle x^{p^{\alpha+\sigma-\gamma}}[y,x]^{p^{\sigma}},
    [y,x]^{p^{\gamma}}, [y,x,x], [y,x,y]\right\rangle,\]
from which we deduce that if $K_3=K_2/N$, then $K_3/Z(K_3)\cong G$, as
desired.
\end{proof}

So we obtain:
\begin{corollary}
Let $G$ be a $2$-generator finite $p$-group of class at most two, $p$
an odd prime. Then $G$ is capable if and only if it is not cyclic and
the orders of the two generators are equal.
\end{corollary}

At this point, an obvious question to ask is whether the necessary
condition of Theorem~\ref{necessity} will also prove to be necessary
for the case of nilpotent groups of class~two, as it was for the
abelian groups. Unfortunately, the answer to that question is no.

Recall that a $p$-group $G$ is ``extra-special'' if and only if
$G'=Z(G)$, $|G'|=p$, and $G^{\rm ab}$ is of exponent~$p$.  A
theorem of Beyl, Felgner, and Schmid in~\cite{beyl} states that an extra-special $p$
group is capable if and only if it is dihedral of order order~$8$, or
of order $p^3$ and exponent $p$, with $p>2$. So the extra-special
$p$-group of order $p^5$ given by
\begin{equation} G = \left\langle x_1,x_2,x_3,x_4 \Biggm| \begin{array}{rcl}
\relax[x_3,x_1][x_3,x_2]^{-1} = [x_3,x_1][x_4,x_1]^{-1} & = & e,\\
\relax[x_4,x_2] = [x_4,x_3] = [x_2,x_1] & = & e,\\
x_1^p=x_2^p=x_3^p=x_4^p & = & e.
\end{array}\right\rangle
\label{eq:extraspecial}
\end{equation}
is not capable, minimally generated by four elements of exponent
$p$. Thus, the necessary condition is not sufficient in general for
groups in~${\germ N}_2$.

A proof that $G$ is not capable, together with some similar
applications using our methods (for example, a proof of Ellis's
Proposition~9 from~\cite{ellis}), will appear in~\cite{capablep}.

For now, we note that the nilpotent product can be used to produce a
natural ``candidate for witness'' to the capability of a given
nilpotent $p$-group $G$. Let $G$ be a finite nilpotent group of
class~$c>0$, minimally generated by elements $x_1,\ldots,x_n$. Let
\[w_1(x_1,\ldots,x_n),\ldots,w_r(x_1,\ldots,x_n)\]
be words in
$x_1,\ldots,x_n$ that give a presentation for $G$, that is:
\[ G = \Bigl\langle x_1,\ldots,x_n \,\Bigm|\,
w_1(x_1,\ldots,x_n),\ldots,w_r(x_1,\ldots,x_n)\Bigr\rangle.\]
Let $\langle y_1\rangle,\ldots,\langle y_n\rangle$ be infinite cyclic
groups, and let 
\[K=\langle y_1\rangle\amalg^{{\germ N}_{c+1}}\cdots \amalg^{{\germ
    N}_{c+1}} \langle y_n\rangle,\]
and let 
\[ N = \Bigl\langle
w_1(y_1,\ldots,y_n),\ldots,w_r(y_1,\ldots,y_n)\Bigr\rangle^K,\]
that is, the normal closure of the subgroup generated by the words
evaluated at $y_1,\ldots,y_n$. Finally, let $M=[N,K]$.

\begin{theorem}
Notation as in the previous paragraph. Then $G$ is capable if and only
if
\[ G \cong (K/M)\bigm/ Z(K/M).\]
\label{th:specialwitness}
\end{theorem}

\begin{proof} We need only prove the ``only if'' part. Note that the
  map that sends $y_1,\ldots,y_n$ to $x_1,\ldots,x_n$ respectively
  gives a well-defined map from $K$ to $G$, and that $K/\langle
  N,K_{c+1}\rangle\cong G$; since $M\subseteq N$, 
\[(K/M)/\langle NM,K_{c+1}M\rangle\cong G,\]
and $\langle NM, K_{c+1}M\rangle$ is central in $K/M$. Therefore,
  $(K/M)/Z(K/M)$ is a quotient of $G$. We need only show that it is
  actually isomorphic to~$G$.

Since $G$ is capable, there is a group $H$ such that $H/Z(H)\cong
G$. Let $h_1,\ldots,h_n$ be elements of $H$ that map to
$x_1,\ldots,x_n$, respectively. Replacing $H$ by $\langle
h_1,\ldots,h_n\rangle$ if necessary, we may assume that
$h_1,\ldots,h_n$ generate $H$. Since $G$ is of class~$c$, $H$ is of
class~$c+1$, and therefore there exists a unique surjective map
from~$K$ (the relatively free $\germ N_{c+1}$ group of rank $n$) to
$H$, mapping $y_1,\ldots,y_n$ to $h_1,\ldots,h_n$, respectively. We
must have that $w_i(h_1,\ldots,h_n)\in Z(H)$ for $i=1,\ldots,r$, so
therefore the map from $K$ to $H$ factors through $K/M$. Since
$Z(K/M)$ maps into $Z(H)$, the induced mapping $K/M\mapsto H\mapsto G$
factors through $(K/M)/Z(K/M)$. Therefore, we have that $(K/M)/Z(K/M)$
has $G$ as a quotient, and is also isomorphic to a quotient of $G$, as
we saw before. The only way this can occur is if $(K/M)/Z(K/M)\cong
G$, as claimed. \end{proof}

In fact, we can do a bit better. By Lemma~2.1 in~\cite{isaacs}, if $G$
is capable then there is a finite group $H$ such that $H/Z(H)\cong
G$. Therefore, choosing $H$ finite above, we can factor the map from
$K$ through a product
\[ \langle y_1\rangle \amalg^{{\germ N}_{c+1}}\cdots \amalg^{{\germ
    N}_{c+1}} \langle y_n\rangle\]
where the groups are now finite cyclic, rather than infinite cyclic. I
believe that they can be chosen so that the order of $y_i$ is the same
as the order of $x_i$ in~$G$; this is the case when the
exponent is $p>c$, but I have not been able to establish
this in the general case.

\section{Related questions\label{sec:relquest}}

By thinking of the direct sum as the first operator in the family of
nilpotent products, we can think of Baer's Theorem for finitely
generated abelian groups as the $k=1$ case of
Theorem~\ref{capabilitynilkprod}. Although, as we noted, it is
possible to take it as a starting point for discussing capability of
other $p$-groups of small class, it seems to me that it is far easier
to use in proving that certain groups are capable than in proving that
certain groups are \textit{not} capable. Even so, some instances of
the latter will appear in~\cite{capablep}, and
Theorem~\ref{th:specialwitness} certainly shows that a careful study
of the nilpotent products may be used to determine that a group is not
capable.

It would be interesting then to see if one can combine the approach
presented here with the ``high-tech'' techniques involving the
epicenter and the nonabelian tensor square of a group. For instance,
every finite $p$-group of class $k$ is a quotient of a $k$-nilpotent
product of cyclic $p$-groups; and it is known that the nonabelian
tensor square of a quotient of $G$ is a quotient of the nonabelian
tensor square of~$G$. What's more, when we mod out by a central
subgroup, it is in general not hard to calculate the kernel of the map
$G\otimes G\to (G/N)\otimes(G/N)$. Perhaps enough information may be
derived on the tensor square of the quotient to determine the
epicenter in such cases? Note that in our applications, we always took
central quotients to reach the groups that witnessed the capability we
were trying to prove. So a somewhat open ended and vague question
is:

\begin{question} Can our approach be combined with the epicenter and
  the nonabelian tensor product to derive new results on capability?
\end{question}

In~\cite{machenry}, T.~MacHenry proved that the cartesian of a
$2$-nilpotent product $A\amalg^{{\germ N}_2} B$ is isomorphic to the
(abelian) tensor product $A^{\rm ab}\otimes B^{\rm ab}$. This
isomorphism proves very useful in investigating questions relating to
amalgams of groups in $\germ N_2$. Since the nonabelian tensor product
is closely connected with commutators (see~\cite{kappecomm}),
it would be useful to know if MacHenry's theorem can be generalised.

\begin{question} Does the nonabelian tensor product give a
  generalisation of MacHenry's Theorem~\cite{machenry}? 
\end{question}

On somewhat more concrete terms, there is plenty of room to expand on
our investigations. If $A$ and $B$ are capable, then so is $A\oplus
B$. Is this also true for the nilpotent products?

\begin{question} Let $G,K\in{\germ N}_k$. If $G$ and~$K$ are capable,
  is $G\amalg^{{\germ N}_k}K$ also capable?
\end{question}

Several works deal with normal forms which can be used to describe
elements of the nilpotent product of cyclic groups; for example,
\cite{struiktwo} deals with the general case of $k=p+1$, and further
work may be found in
\cites{gaglione,waldingerone,wldingeroneadd,waldgag}. We ask,
specifically:

\begin{question} Let $p$ be an odd prime; is it true that a
  $p$-nilpotent product of cyclic $p$-groups of order
  $p^{\alpha_1},\ldots,p^{\alpha_r}$, with $1\leq
  \alpha_1\leq\cdots\leq\alpha_r$, is capable if and only if $r>1$ and
  $\alpha_r\leq \alpha_{(r-1)}+1$?
\end{question}

\begin{question} Is capability of a $k$-nilpotent product of cyclic
  $p$-groups always expressible as a condition of the form
  $\alpha_r\leq \alpha_{(r-1)}+N(k,p)$, for some nonnegative integer
  $N$ which depends on $k$ and $p$?
\end{question}

Although we worked relatively hard to obtain the bound in
Theorem~\ref{necessity}, it is possible that it could be further improved
by a more careful study, similar to the one made
in~\cite{struiktwo}.  There, Struik shows that in the $(p+1)$-nilpotent
product of cyclic $p$-groups, $p$ appears in the denominator of a
term in the multiplication formulas only when the corresponding basic
commutator is either 
\[ [z, \underbrace{y,y,\ldots,y}_{p}\,]\quad\mbox{or}\quad
   [z,y,\underbrace{z,z,\ldots,z}_{p-1}\,].\]
Thus, an Engel condition would remove the difficulties. For example,
in Lemma~\ref{acentralelement}, if we were dealing with a $2$-group of
class~$3$ which also happens to be a $2$-Engel group (so that $[z,y,y]=e$
holds for all $y,z$), then the conclusion will hold with
$\alpha_r=\alpha_{r-1}$. This is easy to verify, since in a $2$-Engel
group any $2$-generated subgroup is of class~$2$. 

So we ask:

\begin{question} For what values of $k$ and $p$ is the bound
  given in Theorem~\ref{necessity} tight?
\end{question}

We know the bound is tight whenever $k<p$, and when $p=2$.

\begin{question} Can the bound found in Theorem~\ref{necessity} be
  improved for groups~$G$ that satisfy an Engel condition? Other
  general conditions?
\end{question}

It is also possible that the inequality may be further strengthened by
using other results in~\cite{easterfield} in which Easterfield also considers
the upper central series to derive bounds.

We have shown that for $2$-generated $p$-groups of class $2$,
$p$ an odd prime, the necessary condition of Theorem~\ref{necessity}
is also sufficient; in~\cite{capablep} we will show this is also the case for
$3$-generated $p$-groups of class~$2$ and exponent~$p$. And we have
exhibited a $4$-generated $p$-group of class~$2$ which satisfies the
necessary condition and is not capable. So we ask:

\begin{question} Is a $3$-generated $p$-group of class two, $p$ an odd
  prime, capable if and only if it satisfies the necessary condition
  of Theorem~\ref{necessity}?
\end{question}

Finally, in relation to Theorem~\ref{th:specialwitness}, we ask:

\begin{question} Suppose $G$ is a finite $p$-group, minimally
  generated by $x_1,\ldots,x_n$, and let $p^{a_i}$ be the order of
  $x_i$. Can we replace $K$ in Theorem~\ref{th:specialwitness} with
  the $(c+1)$-nilpotent product of cyclic groups of order $p^{a_i}$?
  If not, is there some bound $N$, perhaps depending on $c$, $n$, $p$,
  and $a_1,\ldots,a_n$, but not on $G$, such that we can replace $K$
  with the $(c+1)$-nilpotent product of cyclic groups of order $p^N$?
\end{question}

\section*{Acknowledgements}
I became interested in the question of capability after attending a
talk by L.C.~Kappe and M.~Bacon at the AMS Meeting in Baltimore. They
were kind enough to give me a copy of their preprint
of~\cite{baconkappe}, and L.C.~Kappe very patiently answered my
questions over e-mail. I thank them both very much for their time and
help.

In addition, several people were kind enough to help me during the
preparation of this paper. Ed Hook, Robert Israel, and~Joseph
Silverman helped with some of the calculations in the appendix; Bill
Dubuque brought~\cite{granville} to my attention. George Bergman
offered helpful comments.  Avinoam Mann brought~\cite{easterfield} and
several known results to my attention.  Derek Holt, Mike Newman, and
Vasily Bludov from the Group Pub Forum verified some of my
hand-calculations using MAGMA. Finally, Robert Morse used GAP to
verify some of my constructions, answered many of my
questions, and also provided some counterexamples to guesses of mine. I
thank all of them very much for their time and their help.

\appendix
\section{Some results on commutators of powers\label{app:comms}}

In this appendix, we give some results on commutators and commutator
identities needed to establish Theorem~\ref{necessity}. Since the
results may be of independent interest (for example, much less precise
estimates were used in \cite{nildoms}), and the calculations would be
distracting in our main presentation, I decided to place them here at
the end.

To begin with, we recall three consequences of the collection process:

\begin{lemma}[Lemma H1 in~\cite{struikone}]
Let $x,y$ be any elements of a group; let $c_1,c_2,\ldots$ be the
sequence of basic commutators of weight at least two in $x$ and~$[x,y]$,
in ascending order. Then
\begin{equation}
[x^n,y] = [x,y]^nc_1^{f_1(n)} c_2^{f_2(n)}\cdots
c_i^{f_i(n)}\cdots
\label{approxformula}
\end{equation}
where
\begin{equation}
 f_i(n) = a_1\binom{n}{1} + a_2\binom{n}{2} +
\cdots + a_{w_i}\binom{n}{w_i},
\label{formofthefis}
\end{equation}
with $a_i\in\mathbb{Z}$,
and $w_i$ the weight of $c_i$ in $x$ and $[x,y]$.  If the group is
nilpotent, then the expression in {\rm(\ref{approxformula})} gives an
identity, and the sequence of commutators terminates; otherwise,
{\rm(\ref{approxformula})} can be considered as giving a series of
``approximations'' to $[x^n,y]$ modulo successive terms of the lower
central series.
\label{struiklemmah1}
\end{lemma}

\begin{lemma} [Lemma H2 in~\cite{struikone}]
Let $\alpha$ be a fixed integer, and $G$ a nilpotent group of class at
most $n$. If $b_j\in G$ and $r\geq n$, then
\begin{equation}
[b_1,\ldots,b_{i-1},b_i^{\alpha},b_{i+1},\ldots,b_r] =
   [b_1,\ldots,b_r]^{\alpha} c_1^{f_1(\alpha)}
   c_2^{f_2(\alpha)}\cdots
\label{pullingoutexp}
\end{equation}
 where the $c_k$ are commutators in $b_1,\ldots,b_r$ of weight
   strictly greater than $r$, and every $b_j$, $1\leq j\leq r$ appears
   in each commutator $c_k$, the $c_k$ listed in ascending order. The
   $f_i$ are of the form~{\rm (\ref{formofthefis})}, with
   $a_j\in\mathbb{Z}$, and $w_i$ is the weight of $c_i$ (in the $b_i$)
   minus $(r-1)$.
\label{struiklemmah2}
\end{lemma}

We find {\rm(\ref{pullingoutexp})} useful in situations when we have
commutators in some terms, some of which are shown as powers, and we
want to ``pull the exponent out.'' At other times, we will want to
reverse the process and pull the exponents ``into'' a commutator. In
such situations, we use {\rm ({\ref{pullingoutexp})}} to express
$[b_1,\ldots,b_r]^n$ in terms of other commutators. We will call the
resulting identity $(\mathrm{\ref{pullingoutexp}}')$; that is:
\begin{equation}
\relax[b_1,\ldots,b_r]^n =
      [b_1,\ldots,b_{i-1},b_i^n,b_{i+1},\ldots,b_r]\cdots
      v_2^{-f_2(n)}v_1^{-f_1(n)},\tag{$\mathrm{\ref{pullingoutexp}}'$}
\end{equation}
with the understanding that we will only do this in a nilpotent group
so that the formula makes sense.

\begin{lemma}[Theorem 12.3.1 in~\cite{hall}]
Let $x_1,\ldots,x_s$ be any $s$ elements of a group. Let
$c_1,c_2,\ldots$ be the basic commutators in $x_1,\ldots,x_s$ of
weight at least~$2$, written in increasing order. Then
\begin{equation}
(x_1\cdots x_s)^n = x_{1}^n x_{2}^n\cdots x_{s}^n
c_1^{f_1(n)}\cdots c_i^{f_i(n)}\cdots
\label{approxforpower}
\end{equation}
where $f_i(n)$ is of the form~{\rm(\ref{formofthefis})}, with with
$a_j$ integers that depend only on~$c_i$ and not on~$n$, and $w_i$ the
weight of $c_i$ in the $x_j$.  If the group is nilpotent, then
equation {\rm(\ref{approxforpower})} gives an identity of the group,
and the sequence of commutators terminates. Otherwise,
{\rm(\ref{approxforpower})} gives a series of approximations to
$(x_1\cdots x_s)^n$ modulo successive terms of the lower central
series.
\label{struiktheoremh3}
\end{lemma}

The following lemma is easily established by induction on the weight:

\begin{lemma} Let $F(x_1,x_2)$ be the free
  group on two generators. Then every basic commutator of ${\rm
  weight}\geq 3$ is of the form
\begin{equation}
\relax[x_2,x_1,x_1,c_4,\ldots,c_r]\quad\mbox{or}\quad[x_2,x_1,x_2,c_4,\ldots,c_r]
\label{formofbasic}
\end{equation}
where $r\geq 3$, and $c_4,\ldots,c_r$ are basic commutators in $x$
and~$y$ (we interpret $r=3$ to mean that the commutator is of weight
exactly $3$ in $x_1$ and $x_2$).
\label{descriptionbasic}
\end{lemma}

The main idea in our developement is as follows: if we know that
$[z,y^{p^i}]$ centralises $\langle y,z\rangle$ in a group $G\in{\germ
N}_k$, for some prime $p$ and all integers~$i$ greater than or equal
to a given bound~$a$, then we want to prove that a commutator of the
form $[z^{p^n},y]$ is equal to $[z,y^{p^n}]$. To accomplish this, we
observe that a basic commutator of weight $k$ will have exponent
$p^a$, since we may simply use Lemma~\ref{descriptionbasic} and
Proposition~\ref{Widentities}(ii) to pull the exponent into the second
slot of the bracket. An arbitrary commutator of weight $k$ will also
have the same exponent, since $G_k$ is abelian. For a basic commutator
of weight $k-1$ we may use $(\mathrm{\ref{pullingoutexp}}')$ and
Lemma~\ref{struiktheoremh3} and deduce that a sufficiently high power
of $p$ will again yield the trivial element, by bounding below the
power of $p$ that divides the exponents $f_i(p^n)$.  Then we apply
Lemma~\ref{struiktheoremh3} to deal with an arbitrary element of
$G_{k-1}$.  Continuing in this way, we can show that $\langle
y,z\rangle_3$ is of exponent $p^N$ for some large $N$, and we will
obtain the desired result by applying Lemmas~\ref{struiklemmah1}
or~\ref{struiktheoremh3} to $[z^{p^N},y]$ and $[z,y^{p^N}]$. Most of
the work will go into trying to obtain a good estimate on how large
the ``large $N$'' has to be for everything to work.

So our first task is to bound from below the power of $p$ that divides
an expression of the form (\ref{formofthefis}).

We have the following classical result: 

\begin{theorem}[Kummer; see for example~\cite{granville}]
Let $p$ be a prime, and let $n$ and $m$ be positive integers with $n\geq
m$. The exact power of $p$ that divides the binomial coefficient
$\binom{n}{m}$ is given by the number of ``carries'' when we add $n-m$
and $m$ in base~$p$.
\label{kummer}
\end{theorem}

Recall that if $p$ is a prime, and $a$ is a positive integer, 
we let ${\rm ord}_p(a)=n$ if and
only if $n$ is the exact power of $p$ that divides $a$; that is,
$p^n|a$ and $p^{n+1}\not|a$. Formally, we set ${\rm ord}_p(0)=\infty$.
From Kummer's Theorem, we deduce:

\begin{corollary} Let $p$ be a prime, $n$ a positive integer, and $a$
  an integer with $0<a\leq p^n$. Then the exact power of $p$ that
  divides $\binom{p^n}{a}$ is $n-{\rm ord}_p(a)$.
\label{boundbinom}
\end{corollary}

Recall that if $x$ is a real number, then $\lfloor x\rfloor$
  denotes the largest integer smaller than or equal to~$x$.

\begin{corollary} Let $p$ be a prime, $n$ a positive integer, and $m$
  an integer with $0<m\leq p^n$. If $a_1,\ldots,a_m$ are integers,
  then
\begin{equation}
 a_1\binom{p^n}{1} + a_2\binom{p^n}{2} + \cdots +
  a_m\binom{p^n}{m}\label{sumofbinoms}
\end{equation}
is divisible by $p^{n-d}$, where $d$ is the smallest integer such that
$p^{d+1}>m$; that is, $d=\lfloor\log_p(m)\rfloor$.
\label{boundfis}
\end{corollary}

\begin{proof} Write $m=m_0 + m_1p + \cdots + m_dp^d$, with $0\leq
  m_i<p$ and $m_d>0$. Then for all integers $a$ between $1$ and $m$,
  $0\leq {\rm ord}_p(a)\leq d$. Therefore, $n-{\rm ord}_p(a)\geq n-d$,
  so each summand in {\rm(\ref{sumofbinoms})} is divisible by $p^{n-d}$,
  as claimed. 
\end{proof}

\begin{lemma}
Let $k\geq 3$, and let $G\in{\germ N}_k$. Let $p$ be a prime, $a\geq 0$ an
integer, and $y,z\in G$. Suppose that
\begin{equation}
\forall i\geq a,\qquad [z,y^{p^i},y] = [z,y^{p^i},z] = e.\label{identityppower}
\end{equation}
Then $\langle y,z\rangle_k$ is of exponent $p^a$, and $\langle
y,z\rangle_{k-1}$ is of exponent $p^{a+\left\lfloor\frac{1}{p-1}\right\rfloor}$.
\label{inductionbase}
\end{lemma}

\begin{remark} Condition (\ref{identityppower}) holds, for example, if
  $y^{p^a}$ lies in the centralizer of~$z$. This is the situation in
  which we apply our results in Section~\ref{sec:nec}.
\end{remark}

\begin{proof}
If $c$ is a basic commutator of weight exactly $k$ on $y$ and~$z$,
then we may apply $(\mathrm{\ref{pullingoutexp}}')$ to bring the
exponent into the second slot of its expression
as in {\rm (\ref{formofbasic})}, and deduce that $c^{p^a}=e$. Since
$\langle y,z\rangle_k$ is abelian and generated by the basic
commutators, this proves the first part of the statement.

Now let $c\in\langle y,z\rangle_{k-1}$ be a basic commutator of weight
exactly $k-1$. Let $N>0$, and consider $c^{p^N}$. Applying 
$(\mathrm{\ref{pullingoutexp}}')$ to pull the exponent into the $y$ in the
second slot of~$c$ we obtain that
\[ c^{p^N} = \cdots v_2^{-f_2(p^N)} v_1^{-f_1(p^N)}.\]
The only nontrivial $v_i$ are of weight $k$, and so of exponent $p^a$;
and the corresponding exponents are of the form 
\[f_i(p^N)= a_1\binom{p^N}{1} + a_2\binom{p^N}{2}.\]
This expression is divisible by $p^{N-\lfloor \log_p(2)\rfloor}$, and
so will be trivial whenever $N\geq a+\lfloor \log_p(2)\rfloor$. Therefore,
every basic commutator of weight exactly $k-1$ is of exponent
$p^{a+\lfloor \log_p(2)\rfloor}$. Since $k>2$, $G_{k-1}$ is abelian as
well; every generator is of exponent $p^{a+\lfloor \log_p(2)\rfloor}$;
the theorem is now proven by noting that $\lfloor\log_p(2)\rfloor =
\left\lfloor\frac{1}{p-1}\right\rfloor$. 
\end{proof}

We will need the following small calculation:

\begin{lemma}
Let $k,n$ be positive integers, $n>1$. The maximum of the numbers
\[ \left\lfloor\frac{k-s}{n-1}\right\rfloor +
\lfloor\log_n(s+1)\rfloor,\quad s=1,\ldots,k,\]
is equal to $\left\lfloor\frac{k}{n-1}\right\rfloor$, and it is always
attained at $s=n-1$.
\label{boundfors}
\end{lemma}

\begin{proof} If $k<n-1$, the result is immediate. Assume $k\geq n-1$, and write $k=q(n-1)+r$, with $0\leq
  r<n-1$. If $r>0$, then the first summand is equal to $q$ for the
  first $r$ terms, while the second summand is~$0$. Then the first
  summand drops by $1$ and the second summand stays at zero until we
  get to $s=n-1$, when it becomes equal to~$1$; after that, the first
  summand drops faster than the second summand increases. Thus, the
  maximum is indeed obtained at $s=n-1$ (as well as at all values
  $s\leq r$ and a few after $s=n-1$). If $r=0$, then for
  $s=1,\ldots,n-1$ the first summand is equal to $q-1$, but only at
  $s=n-1$ is the second summand, $\lfloor \log_n(s)\rfloor$,
  positive. So in this case the maximum is attained at $s=p-1$, and
  only there.  In either case, the value of that maximum is:
\[ \left\lfloor \frac{k-(n-1)}{n-1}\right\rfloor +
  \lfloor\log_n(n)\rfloor = \left\lfloor \frac{k}{n-1}-1\right\rfloor
  + 1 = \left\lfloor\frac{k}{n-1}\right\rfloor,\]
as claimed.
\end{proof}

\begin{lemma}[cf. Lemma~8.83 in~\cite{nildoms}] 
Let $k>0$, and let $G\in{\germ N}_k$ be a group. Let $p$ be a prime,
$a>0$ an integer, and $y,z\in G$. Assume that $a$, $y$, $z$, and $G$
satisfy {\rm (\ref{identityppower})}.
Then $\langle y,z\rangle_{k-m}$ is of exponent $p^{a+\left\lfloor\frac{m}{p-1}\right\rfloor}$
for $m=0,1,\ldots,(k-3)$.
\label{expofthirdterm}
\end{lemma}

\begin{proof} The result is vacuously true if $k=1$ or $k=2$.  We proceed by induction on
  $m$. Lemma~\ref{inductionbase} proves cases $m=0$ and $m=1$, so we also assume $m\geq
  2$.

Assume the result is true for $n=0,1,\ldots,m-1$. First, we consider a
basic commutator $c$ of weight exactly $k-m$. We may express $c$ as in
{\rm (\ref{formofbasic})}, and applying
$(\mathrm{\ref{pullingoutexp}}')$ to $c^{p^N}$ we obtain
\[ c^{p^N} = [z,y^{p^N},w,c_4,\ldots,c_r]\cdots v_2^{-f_2(p^N)}
v_1^{-f_1(p^N)},\]
where $w\in\{y,z\}$, $c_4,\ldots,c_r$ are basic commutators in $z$ and
$y$, and the $v_i$ are basic commutators in $z,y,w,c_4,\ldots,c_r$;
we know that each of $z,y,w,c_4,\ldots,c_r$ appears at least once in
each $v_i$, and that we may express the weight of $v_i$ on
$z,y,w,c_4,\ldots,c_r$
as $r+s$ for some positive integer~$s$. Thus,
we may conclude that if the weight of $v_i$ in $z,y,w,c_4,\ldots,c_r$
is $r+s$, then $W(v_i)\geq (k-m)+s$. In particular, we consider only those
commutators with $1\leq s\leq m$. 

If $v_i$ is of weight $r+s$ in $z,y,w,c_4,\ldots,c_r$, then by
Lemma~\ref{struiklemmah2} we have that
\[ f_i(p^N) = a_1\binom{p^N}{1}+\cdots +a_{s+1}\binom{p^N}{s+1}.\]
Therefore, $f_i(p^N)$ is divisible by $p^{N-\lfloor
\log_p(s+1)\rfloor}$. Since $v_i\in\langle y,z\rangle_{(k-(m-s))}$, by
the induction hypothesis we know that $v_i^{-f_i(p^N)}$ is trivial whenever
$N-\lfloor\log_p(s+1)\rfloor \geq a+\left\lfloor\frac{m-s}{p-1}\right\rfloor$. Therefore, if
\[ N \geq a + \left\lfloor\frac{m-s}{p-1}\right\rfloor + \lfloor\log_p(s+1)\rfloor\]
for $s=1\ldots,m$, then we may conclude that $c^{p^N}=e$. By
Lemma~\ref{boundfors}, the greatest of these values is
$a+\lfloor\frac{m}{p-1}\rfloor$, which shows that
the basic commutators of weight exactly $m-k$ are of
exponent $p^{a+\lfloor\frac{m}{p-1}\rfloor}$, as desired.

Now take an arbitrary element of $\langle y,z\rangle_{(k-m)}$, and
write $c = d_1\cdots d_r$, where $d_i=c_i^{a_i}$ is the power of a
basic commutator of weight at least $k-m$ in $y$ and~$z$, and
$c_1<\cdots< c_r$. To estimate $c^{p^N}$ we apply
Lemma~\ref{struiktheoremh3} to this expression. We obtain:
\[(d_1\cdots d_r)^{p^N} = c_{1}^{p^N} c_{2}^{p^N}\cdots c_{r}^{p^N}
u_1^{f_1(p^N)}\cdots u_i^{f_i(p^N)}\cdots
\]
where $f_i(p^N)$ is of the form~{\rm(\ref{formofthefis})}, with $w_i$
the weight of $u_i$ in the $d_j$. If we let the weight of $u_i$ in the
$d_j$ be equal to $s$, then we know that $2\leq s$.  Since $W(d_i)\geq
(m-k)$, we have that $W(u_i)\geq s(m-k)$. In particular, we may
restrict to values of $s$ satisfying $2\leq s\leq
\left\lfloor\frac{k}{k-m}\right\rfloor$.  

If $u_i$ is of weight $s$ in the $d_j$, then it lies in $\langle
y,z\rangle_{k-(m-(s-1)(k-m))}$, so by the induction hypothesis it is
of exponent $p^{a+\lfloor\frac{m-(s-1)(k-m)}{p-1}\rfloor}$. And we know that
\[ f_i(p^N) = a_1\binom{p^N}{1} + \cdots + a_s\binom{p^N}{s}.\]
This is a multiple of $p^{N-\lfloor\log_p(s)\rfloor}$, so we can
guarantee that $u_i^{f_i(p^N)}$ is trivial if
\[ N \geq a + \left\lfloor\frac{m-(s-1)(k-m)}{p-1}\right\rfloor + \lfloor\log_p(s)\rfloor.\]
Since $s\geq 2$ and $k-m\geq 3$, it is clear that this value will
certainly be no larger than $a+\lfloor\frac{m}{p-1}\rfloor$, so we may
conclude that any element of $\langle y,z\rangle_{(k-m)}$ is of
exponent $p^{a+\lfloor\frac{m}{p-1}\rfloor}$, as claimed. This
finishes the induction.
\end{proof}

\begin{theorem}[cf. Corollary~8.84 in~\cite{nildoms}]
Let $k\geq 1$ and $G\in{\germ N}_{k+1}$. Let $p$ be a prime, $a>0$ an
integer, and $y$ and~$z$ elements of~$G$. Assume that
\[ \forall i\geq a,\qquad [z,y^{p^i},y]=[z,y^{p^i},z]=e.\]
If $N\geq a+\left\lfloor\frac{k-1}{p-1}\right\rfloor$, then
$[z^{p^N},y] = [z,y]^{p^N} = [z,y^{p^N}]$.
\label{likenil2withpowers}
\end{theorem}

\begin{proof} We apply Lemma~\ref{struiklemmah1} to $[z^{p^N},y]$
  (applying Lemma~\ref{struiktheoremh3} yields the same result):
\[ [z^{p^N},y] = [z,y]^{p^N} v_1^{f_1(p^N)}v_2^{f_2(p^N)}\cdots,\]
where $v_1,v_2,\ldots$ are the basic commutators of weight at least
$2$ in $z$ and $[z,y]$. If $v_i$ is of weight $s\geq 2$ in $z$ and
$[z,y]$, then 
\[ f_i(p^N) = a_1\binom{p^N}{1} + a_2\binom{p^N}{2} + \cdots +
a_s\binom{p^N}{s},\]
and we know that $v_i\in\langle y,z\rangle_{s+1}$. By
Theorem~\ref{expofthirdterm} $v_i^{f_i(p^N)}$ will be trivial if
\[ N \geq a+\left\lfloor\frac{(k+1)-(s+1)}{p-1}\right\rfloor + \lfloor\log_p(s)\rfloor.\]
This must hold for $s=2,\ldots,k$. We set $t=s-1$ and rewrite it as:
\[ N \geq a + \left\lfloor\frac{(k-1)-t}{p-1}\right\rfloor +
\lfloor\log_p(t+1)\rfloor,\]
with $t=1,\ldots,k-1$
By Lemma~\ref{boundfors}, the largest of these numbers is 
\[ a + \left\lfloor\frac{k-1}{p-1}\right\rfloor,\]
which proves that if $N\geq
a+\left\lfloor\frac{k-1}{p-1}\right\rfloor$, then
$[z^{p^N},y]=[z,y]^{p^N}$. The proof that $[z,y^{p^N}]=[z,y]^{p^N}$ is
essentially the same.
\end{proof}

\begin{remark} The results above also hold if the identities
  in~$(\mathrm{\ref{identityppower}})$ are replaced by an identity in
  which $p^i$ is placed on any of the three slots in $[z,y,z]$, and in
  any of the three slots in $[z,y,y]$; for instance,
\[ \relax[z^{p^i},y,z]=[z,y,y^{p^i}]=e,\quad\mbox{or}\quad
  [z,y^{p^i},z]=[z^{p^i},y,y]=e,\]
etc.
\end{remark}

\end{document}